\documentclass{amsproc}

\author{Frances Kirwan}
\address{Mathematical Institute, Oxford University, Oxford OX1 3LB, UK}
\email{kirwan@maths.ox.ac.uk}

\thanks{The author is a member of VBAC (Vector Bundles on Algebraic Curves),
which is partially supported by EAGER (EC FP5 Contract no. HPRN-CT-2000-00099),
and acknowledges with gratitude the hospitality of the University of Melbourne 
and the Australian National University during
the writing of part of this paper.}



\usepackage{latexsym}

\usepackage{euscript}

\title[The Yang--Mills Stratification Revisited]{Moduli Spaces of Bundles over Riemann Surfaces and the Yang--Mills Stratification Revisited}

\newtheorem{prop}{Proposition}[section]
\newtheorem{lem}[prop]{Lemma}
\newtheorem{cor}[prop]{Corollary}
\newtheorem{thm}[prop]{Theorem}

\newtheorem{dff}[prop]{Definition}
\newenvironment{df}{\begin{dff}\rm}{\end{dff}}
\newtheorem{REM}[prop]{Remark}
\newenvironment{rem}{\begin{REM}\rm}{\end{REM}}
\newtheorem{examplit}[prop]{Example}
\newenvironment{example}{\begin{examplit}\rm}{\end{examplit}}

\newcommand{\mnd}{\mathcal{M}(n,d)}
\newcommand{\Mnd}{\mathcal{M}(n,d)}

\newcommand{\HS}{H^{*}}
\newcommand{\ar}{a_{r}}
\newcommand{\bjr}{b_{r}^{j}}
\newcommand{\fr}{f_{r}}

\newcommand{\V}{\hat{V}}
\newcommand{\G}{\mathcal{G}}
\newcommand{\GG}{\overline{\mathcal{G}}}

\newcommand{\Cmu}{\mathcal{C}_{\mu}}
\newcommand{\HG}{H^{*}_\mathcal{G}}
\newcommand{\Css}{\mathcal{C}^{ss}}
\newcommand{\C}{\mathcal{C}}
\newcommand{\liet}{{\bf t}}

\newcommand{\stab}{{\rm Stab}}

\newcommand{\RR}{{\Bbb R }}
\newcommand{\CC}{{\Bbb C }}
\newcommand{\ZZ}{{\Bbb Z}}
\newcommand{\PP}{ {\Bbb P } }
\newcommand{\QQ}{{\Bbb Q }}

\newcommand{\renorm}{{ \setcounter{equation}{0} }}

\begin{document}

\begin{abstract}
Refinements of the Yang--Mills stratifications of spaces of connections
over a compact Riemann surface $\Sigma$ are investigated.
The motivation for this study is the search for a complete set of relations
between the standard generators for the cohomology of the moduli
spaces $\mnd$ of stable holomorphic bundles of rank $n$ and degree
$d$ when $n$ and $d$ are coprime and $n>2$.
\end{abstract}

\maketitle

\setcounter{page}{1}

\newcommand{\n}{\hat{n}}
\newcommand{\D}{\hat{d}}
\newcommand{\M}{\mathcal{M}(\n,\D)}

The moduli space
$\mnd$ of semistable holomorphic bundles of rank $n$ and degree $d$ over a Riemann
surface $\Sigma$ of genus $g \geq 2$ can be constructed as a quotient of an
infinite dimensional affine space of connections $\mathcal{C}(n,d)$ by a complexified gauge
group $\G_c(n,d)$,
in an infinite-dimensional version of the 
construction of quotients in Mumford's geometric invariant theory \cite{MFK}.
When $n$ and $d$ are coprime, $\mnd$ is the topological quotient of the semistable
subset $\mathcal{C}(n,d)^{ss}$ of $\mathcal{C}(n,d)$ by the action of $\G_c(n,d)$.
Any nonsingular complex projective variety on which a complex reductive
group $G$ acts linearly has a $G$-equivariantly perfect stratification
by locally closed nonsingular $G$-invariant subvarieties with its set of semistable
points $X^{ss}$ as an open stratum. This stratification can be obtained as the
Morse stratification for the normsquare of a moment map on $X$ \cite{K2}; in the 
case of  the moduli space $\Mnd$ 
  the r\^{o}le of
the normsquare of the moment map is played by the Yang-Mills
functional. In \cite{K7}
this Morse stratification of $X$ is refined to obtain a stratification of $X$ by locally closed
nonsingular $G$-invariant subvarieties with the set $X^s$ of stable points of $X$ as an open stratum. 
The other strata can be defined inductively in terms of the sets of stable
points of closed nonsingular subvarieties of $X$, acted on by reductive
subgroups of $G$, and their projectivised normal bundles.

In their fundamental paper \cite{AB}, Atiyah and Bott studied
a stratification of $\mathcal{C}(n,d)$ defined using the Harder-Narasimhan
type of a holomorphic bundle over $\Sigma$, which they expected 
to be the Morse stratification of the Yang-Mills 
functional (this was later shown to be the case \cite{Daskal}).
The aim of this paper is to apply the methods of \cite{K7} to the 
Yang-Mills stratification to obtain refined stratifications of $\mathcal{C}(n,d)$,
and to relate these stratifications to natural refinements of the notion of the Harder-Narasimhan
type of a holomorphic
bundle over $\Sigma$. 

The
motivation for this study was the search for a complete set of relations
among the standard generators for the cohomology of the moduli
spaces $\Mnd$ when the rank $n$ and degree $d$ are coprime and $n>2$ \cite{EK}.
The cohomology rings of the moduli spaces $\mnd$  have been the subject of much interest
over many years; see for example \cite{AB,Be,DR,D,G,HN,HS,JK2,KN,LP,NR,NS,T,WI,Z}
among many other pieces of work. In the case when $n=2$ we now have
 a very thorough understanding of the
structure of the cohomology ring
\cite{B,KN,ST,Z}. For arbitrary $n$ it is known that the cohomology has no
torsion and formulas for computing the Betti
numbers have been obtained,
 as well as a set of
generators for the cohomology ring  \cite{AB,DB,DR,HN,Z2}. 
 When $n=2$ the
relations between these generators can be explicitly described 
and in particular a conjecture of Mumford, that a
certain set of relations is a complete set, was proved some years ago
\cite{B,KN,K,ST,Z}.
However less is known about the relations between the
generators when $n>2$, and the most
obvious generalisation of Mumford's conjecture to the cases when $n>2$
is false,
 although a modified version of the conjecture (using `dual
Mumford relations' together with the original Mumford relations) is
true for $n=3$ \cite{E}. 
There is, however, a further generalisation of
Mumford's relations,
and the attempt by Richard Earl and the author to prove that this set is indeed complete 
was the original stimulus for studying refinements of the Yang-Mills
stratification, although a different application has now appeared \cite{EK,JKKW2}. 

The layout of this paper is as follows. $\S$1 recalls background
material on moduli spaces of bundles and different versions
of Mumford's conjecture, while $\S$2 reviews the Morse stratification
of the normsquare of a moment map and some refinements
of this stratification.  $\S$3 studies the structure of subbundles of semistable bundles
over $\Sigma$ which
are direct sums of stable bundles all of the same slope, and this is
used in $\S$4 to define two canonical refinements of the Harder-Narasimhan
filtration of a holomorphic bundle over $\Sigma$,  and thus
to construct two refinements of the Yang-Mills stratification. In the next two
sections the stratification defined in $\S$2 is applied to holomorphic
bundles over $\Sigma$; its indexing set is studied in $\S$5 and
the associated strata are investigated in $\S$6. This stratification corresponds
to a third refinement of the Harder-Narasimhan filtration whose subquotients
are all direct sums of stable bundles of the same slope. The relationship
between these three filtrations is considered in $\S$7, and 
$\S$8 provides a brief conclusion.

\section{Background material on moduli spaces of bundles}
\renorm

When $n$ and $d$ are coprime, 
the generators for the rational cohomology\footnote{In 
this paper all cohomology will have rational coefficients.}
  of the moduli space
$\mnd$ given by Atiyah and Bott in 
\cite{AB} are
obtained from a (normalised) universal bundle $V$ over $\mnd \times
\Sigma$. With respect to the K\"{u}nneth decomposition of 
$
\HS(\mnd \times \Sigma)
$
the $r$th Chern class $c_{r}(V)$ of $V$ can be
written as
$$
c_{r}(V) = a_{r} \otimes 1 + \sum_{j=1}^{2g} b_{r}^{j} \otimes
\alpha_{j} + f_{r} \otimes \omega
$$
where $\{1\}, \{\alpha_{j}:1 \leq j \leq 2g\},$ and $\{ \omega\}$ are
standard bases for $H^0(\Sigma),H^1(\Sigma)$ and $H^2(\Sigma)$, and
\begin{equation}
a_{r} \in H^{2r}(\mnd), \qquad b_{r}^{j} \in H^{2r-1}(\mnd), \qquad
f_{r} \in H^{2r-2}(\mnd), \label{gen}
\end{equation}
for $1 \leq r \leq n$ and $1\leq j \leq 2g$. It was shown by Atiyah and
Bott 
\cite[Prop. 2.20 and p.580]{AB} 
that the classes $a_{r}$ and $f_{r}$ (for $2 \leq r \leq
n$) and $b_{j}^r$ (for $1 \leq r \leq n$ and $1 \leq j \leq 2g$)
generate the rational cohomology ring of $\mnd$.\\
\indent Since tensoring by a fixed holomorphic line bundle of degree
$e$ gives an isomorphism between the moduli spaces $\mnd$ and $\mathcal{
M}(n,d+ne)$, we may assume without loss of generality that
$$
(2g-2)n < d < (2g-1)n.
$$
This implies that $H^1(\Sigma ,E) = 0$ for any stable bundle of rank
$n$ and degree $d$ \cite[Lemma 5.2]{N2},
 and hence that $\pi_{!}V$ is
a bundle of rank $d-n(g-1)$ over $\mnd$, where 
$$
\pi: \mnd \times \Sigma \to \mnd 
$$
is the projection onto the first component and $\pi_!$ is the
K-theoretic direct image map. It follows that 
$$
c_{r} (\pi_! V) = 0
$$
for $r>d-n(g-1)$. Via the Grothendieck-Riemann-Roch theorem we can
express the Chern classes of $\pi_! V$ as polynomials in the
generators $a_r, b_r^j, f_r$ described above, and hence their
vanishing gives us relations between these generators.  These
are Mumford's relations, and they give us a complete set of relations
when $n=2$. We can generalise them for $n>2$ as follows.

\indent Suppose that $0<\n<n$, and that $\D$ is coprime to $\n$. Then
we have a universal bundle $\V$ over $\M \times \Sigma$, and both $V$
and $\V$ can be pulled back to $\M \times \mnd \times \Sigma$. If
$$
\frac{\D}{\n} > \frac{d}{n}
$$
then there are no non-zero holomorphic bundle maps from a stable
bundle of rank $\n$ and degree $\D$ to a stable bundle of rank $n$ and
degree $d$, and hence, if 
$$
\pi: \M \times \mnd \times \Sigma \to \M \times \mnd,
$$
is the projection onto the first two components, it follows that 
$
-\pi_! (\V^* \otimes V)
$
is a bundle of rank $n\n(g-1)-d\n+\D n$ over $\M \times \mnd$. Thus 
$$
0 = c_{r} (-\pi_! (\V^* \otimes V)) \in  \HS(\M \times \mnd)
$$
if $r>n\n(g-1)-d\n+\D n$ and hence the slant product
$$
c_{r} (-\pi_! (\V^* \otimes V)) \backslash \gamma \in \HS(\mnd)
$$
of $c_{r} (-\pi_! (\V^* \otimes V))$ with any homology class $\gamma
\in H_*(\M)$ vanishes when 
$$
r>n\n(g-1)-d\n+\D n.
$$
The relations 
between the generators $a_r,b_r^j,f_r$ obtained in this way for $0<\n<n$ and 
$$
\frac{d}{n}+1 > \frac{\D}{\n} > \frac{d}{n}
$$
and 
$$
n\n(g-1)-d\n+\D n <r < n\n(g+1)-d\n+\D n
$$
(with a little more care taken when $\n$ and $\D$ are not coprime)
are the ones we consider. They 
are essentially Mumford's
relations when $n=2$.
To show that these form a complete set of
relations, a natural strategy is to consider the Yang-Mills stratification which
was used by Atiyah and Bott to obtain their generators for the cohomology ring.

Recall that a holomorphic vector bundle $E$ over $\Sigma$ is
called {\em  semistable} (respectively {\em stable}) if every holomorphic
subbundle $D$ of $E$ satisfies
$$
\mu (D) \leq \mu(E) , \indent (\mbox{respectively } \mu(D) < \mu(E)),
$$
where $\mu(D) =$ degree($D$)/rank($D$) is the {\em slope} of
$D$. Bundles which are not semistable are said to be {\em unstable}. Note that
semistable bundles of coprime rank and degree are stable.

Let $\mathcal{E}$ be a fixed $C^{\infty}$ complex vector bundle
of rank $n$ and degree $d$ over $\Sigma$. Let $\C$ be the space of all holomorphic
structures on $\mathcal{E}$ and let $\G_{c}$ denote the group of all
$C^{\infty}$ complex automorphisms of $\mathcal{E}$. Atiyah and Bott
\cite{AB}
 identify the moduli space
$\mnd$ with the quotient $\Css/\mathcal{G}_{c}$ where $\Css$ is the open
subset of $\C$ consisting of all semistable holomorphic structures on
$\mathcal{E}$. The group $\G_{c}$ is the
complexification of the gauge group $\G$ which consists of all smooth automorphisms
of $\mathcal{E}$ which are unitary with respect to a fixed Hermitian
structure on $\mathcal{E}$. We shall write
$\overline{\G}$ for the quotient of $\G$ by its $U(1)$-centre and
$\overline{\G}_{c}$ for the quotient of $\G_{c}$ by its ${\bf 
C}^{*}$-centre. There are natural isomorphisms
$$
H^*(\mnd) \cong H^{*}(\Css/\G_{c}) = H^{*}(\Css/\overline{\G}_{c}) \cong
H^{*}_{\overline{\G}_{c}}(\Css) \cong H^{*}_{\overline{\G}}(\Css)
$$
between the cohomology of the moduli space and the 
$\overline{\G}$-equivariant cohomology of $\Css$, 
since the ${\bf  C}^{*}$-centre of $\G_{c}$ acts trivially on $\Css$, while 
$\overline{\G}_{c}$ acts freely on $\Css$ and $\overline{\G}_{c}$ is
the complexification of $\overline{\G}$. 

In order to show that the restriction map $H^{*}_{\overline{\G}}(\C)
\rightarrow H^{*}_{\overline{\G}}(\Css)$ is surjective, Atiyah
and Bott consider the Yang--Mills (or Atiyah--Bott--Shatz) stratification
of $\C$.
This stratification $\{\Cmu : \mu \in \mathcal{M} \}$ is
the Morse stratification for the Yang--Mills functional on
$\C$, but it also has a more explicit description. It is
indexed by the partially ordered set $\mathcal{M}$ consisting of all the Harder--Narasimhan
{\em types} of holomorphic bundles of rank $n$ and degree $d$, defined as follows.
Any holomorphic bundle $E$ over $M$ of rank $n$ and degree $d$ has a
filtration 
$$
0 = E_{0} \subset E_{1} \subset \cdot \cdot \cdot \subset E_{P} = E
$$
of subbundles such that the subquotients $Q_{p}= E_{p}/E_{p-1}$
are semi-stable for $1 \leq p \leq P$ and satisfy 
$$
\mu(Q_{p}) = \frac{d_{p}}{n_{p}} > \frac{d_{p+1}}{n_{p+1}} =
\mu(Q_{p+1})
$$ 
where $d_{p}$ and $n_{p}$ are the degree and rank of
$Q_{p}$ and $\mu(Q_p)$ is its slope. This filtration is canonically associated to $E$ and is called 
the Harder-Narasimhan filtration of $E$. We define the {\em type} of $E$ to be
$$
\mu = (\mu(Q_{1}),...,\mu(Q_{P})) \in \QQ^{n}
$$
where the entry $\mu(Q_{p})$ is repeated $n_{p}$ times. The
semistable bundles have type 
$ \mu_{0} = (d/n,...,d/n) $
and form the
unique open stratum. The set $\mathcal{M}$ of all possible types of holomorphic
vector bundles over $\Sigma$ provides our indexing set, and if 
$\mu \in \mathcal{M}$ the subset $\Cmu \subseteq \mathcal{C}$ is
defined to be the set of all holomorphic vector bundles over $\Sigma$ of type
$\mu$. 
A partial order on
$\mathcal{M}$ with the property that the closure of the stratum indexed by $\mu$ is
contained in the union of all strata indexed by $\mu' \geq \mu$ is defined as follows. 
Let $\sigma=(\sigma_{1},...,\sigma_{n})$ and
$\tau=(\tau_{1},...,\tau_{n})$ be two types; then
\begin{equation} \label{po}
\sigma \geq \tau \mbox{ if and only if }
\sum_{j \leq i} \sigma_{j} \geq \sum_{j \leq i} \tau_{j} \mbox{ for } 1 \leq i
\leq n-1.
\end{equation}
The gauge group $\G$ acts on $\C$ preserving the
stratification which is equivariantly
perfect with respect to this action, which means that its equivariant Thom-Gysin sequences
$$
\cdots \to H_{\G}^{j-2d_\mu}(\Cmu) \to H^j_{\G}(U_\mu) \to
H^j_{\G}(U_\mu-\Cmu) \to \cdots 
$$ break up into short
exact sequences
$$
0 \to H_{\G}^{j-2d_\mu}(\Cmu) \to H^j_{\G}(U_\mu) \to
H^j_{\G}(U_\mu-\Cmu) \to 0.
$$
Here 
\begin{equation}
d_{\mu}= \sum_{i>j} (n_{i}d_{j}-n_{j}d_{i}+n_{i}n_{j}(g-1)), \label{12}
\end{equation}
is the complex codimension of $\Cmu$ in $\C$ and 
$U_\mu$ is the open subset of $\C$ which is the union of all
those strata labelled by $\mu'\leq \mu$; we also have
$$ H^*_\G (\Cmu) \cong \bigotimes_{j=1}^P H^*_{\G(n_j,d_j)}(\C(n_j,d_j)^{ss}). 
$$
Atiyah and Bott show that the stratification is equivariantly perfect
by considering the composition   
of the
Thom-Gysin map
$H_{\G}^{j-2d_\mu}(\Cmu) \to H^j_{\G}(U_\mu)$
with restriction to $\Cmu$, 
which is multiplication by the equivariant Euler class $e_\mu$ of the
normal bundle to $\Cmu$ in $\C$. They show 
that $e_\mu$ is not a zero-divisor in $\HG(\Cmu)$ and deduce that the
Thom-Gysin maps 
$H_{\G}^{j-2d_\mu}(\Cmu) \to H^j_{\G}(U_\mu)$
are all injective.

So putting this all together Atiyah and Bott obtain inductive formulas for the $\G$-equivariant
Betti numbers of $\C^{ss}$, and they also conclude that there is a natural surjection
\begin{equation}
H^{*}(B\overline{\G}) \cong H^{*}_{\overline{\G}}(\C)
\rightarrow H^{*}_{\overline{\G}}(\Css) \cong H^{*}(\mnd). \label{14}
\end{equation}
Thus generators of the cohomology ring
$\HS(B\overline{\G})$ give generators of the cohomology ring $\mnd$.

The classifying space $B\G$
can be identified with the space $\mbox{Map}_{d}(\Sigma,BU(n))$ of all
smooth maps $f:\Sigma \rightarrow BU(n)$ such that the pullback to $\Sigma$ of
the universal vector bundle over $BU(n)$ has degree $d$. If we
pull back this universal bundle using the evaluation map
$$
\mbox{Map}_{d}(\Sigma,BU(n)) \times \Sigma \rightarrow BU(n): (f,m) \mapsto f(m)
$$
then we obtain a rank $n$ vector bundle $\mathcal{V}$ over $B\G \times
\Sigma$. If further we restrict the pullback bundle induced by the maps
$$
\Css \times E\G \times \Sigma \rightarrow \C \times E\G \times \Sigma \rightarrow \C
\times_{\G} E\G \times \Sigma \stackrel{\simeq}{\rightarrow} B\G \times \Sigma
$$
to $\Css \times \{e\} \times \Sigma$ for some $e \in E\G$ then we obtain
a $\G$-equivariant holomorphic bundle on $\Css \times \Sigma$. The
${\bf C}^*$-centre of
$\G$ acts as scalar multiplication on the fibres,
and the associated projective bundle descends to a holomorphic
projective bundle over $\mnd \times \Sigma$.
 In fact this is the projective bundle of a
holomorphic vector bundle $V$ over $\mnd \times \Sigma$
which has the universal
property that, for any $[E] \in \mnd$ representing a bundle $E$ over
$\Sigma$, the restriction of $V$ to $\{[E]\} \times \Sigma$ is
isomorphic to $E$.\\
\indent By a slight abuse of notation we define elements $\ar, \bjr, \fr$
in $\HS(B\G;\QQ)$ by writing
$$
c_{r}(\mathcal{V}) = \ar \otimes 1 + \sum_{j=1}^{2g} \bjr \otimes
\alpha_{j} + \fr \otimes \omega \indent 1 \leq r \leq n.
$$
where, as before, $\omega$ is the standard generator of $H^{2}(\Sigma)$ and
$\alpha_{1},...,\alpha_{2g}$ form a fixed canonical cohomology basis for
$H^{1}(\Sigma)$.
In fact the ring
$H^{*}(B\G)$ is freely generated 
as a graded super-commutative algebra over $\QQ$ by the elements 
$$
\{\ar : 1 \leq r \leq n\} \cup \{\bjr : 1 \leq r \leq n, 1 \leq j \leq 2g\} \cup 
\{ \fr : 2 \leq r \leq n\}
$$
and if we omit $a_1$ we get  $H^{*}(B\overline{\G})$. These generators
restrict to the generators $\ar, \bjr, \fr$ given at (\ref{gen}) for
$\HS(\mnd)$ under the surjection (\ref{14}).\\
\indent The relations among these generators for $\HS(\mnd;\QQ)$ are
then given by the kernel of the restriction map (\ref{14}) which is in
turn determined by the map
\begin{equation}
\HG(\C) \cong  H^{*}_{\overline{\G}}(\C) \otimes \HS(BU(1)) \rightarrow
H^{*}_{\overline{\G}}(\Css) \otimes
\HS(BU(1)) \cong \HG(\Css), \label{surj}
\end{equation}
and the proof that the Yang--Mills stratification is equivariantly perfect leads to 
{ completeness criteria} for a set of relations to be
complete.
Let $\mathcal{R}$ be a
subset of the kernel of the restriction map 
$\HG(\C) \rightarrow \HG(\Css).$
Suppose that for each unstable type $\mu \neq \mu_{0}$ there is a subset 
$ \mathcal{R}_{\mu}$ of the ideal generated by $\mathcal{R}$ in $\HG(\C)$ such that
the image of $\mathcal{R}_{\mu}$ under the restriction map
$$
\HG(\C) \rightarrow \HG(\C_{\nu})
$$
is zero unless $\nu \geq \mu$ and when $\nu = \mu$ contains the ideal of
$\HG(\C_{\mu})$ generated by the equivariant
Euler class $e_{\mu}$ of the normal bundle to the stratum $\C_{\mu}$ in
$\C.$ Then $\mathcal{R}$ generates the kernel of the restriction map
$
\HG(\C) \rightarrow \HG(\Css)
$
as an ideal of $\HG(\C).$

In fact
Atiyah and Bott could have replaced the Yang-Mills stratification with a coarser stratification of
$\C$ and obtained equivalent results.
For any integers $n_1$ and $d_1$ let $S_{n_1,d_1}$ be the
subset of $\C$ consisting of all those holomorphic structures with
Harder-Narasimhan filtration 
$
0 = E_0 \subset E_1 \subset \cdots \subset E_s = E
$
where $E_1$ has rank $n_1$ and degree $d_1$. We shall say that such a
holomorphic structure has {\em coarse type} $(n_1,d_1)$.
Then $S_{n_1,d_1}$ is
locally a submanifold of finite codimension 
\begin{eqnarray*}
\delta_{n_1,d_1} & = & n d_1 - n_1 d + n_1 (n-n_1)(g-1)
\end{eqnarray*}
in $\C$
and
\begin{equation}
\label{tensor}
\quad \HG(S_{n_1,d_1}) \cong \HS_{\G(n_1,d_1)}(\C(n_1,d_1)^{ss}) \otimes 
 \HS_{\G(n-n_1,d-d_1)}\left( U(n_1,d_1) 
 \right)
\end{equation}
where 
$$ U(n_1,d_1) =
\bigcup_{\frac{d_2}{n_2} < \frac{d_1}{n_1}} S(n-n_1,d-d_1)_{n_2,d_2}
$$
is an open subset of $\C(n-n_1,d-d_1)$. Moreover the equivariant Euler class $e_{n_1,d_1}$ of
the normal to $S_{n_1,d_1}$ in $\C$ is not a zero divisor in
$\HG(S_{n_1,d_1})$, so the stratification of $\C$ by coarse type
$$
\left\{ S_{n_1,d_1}: 0 < n_1 < n, \frac{d_1}{n_1} > \frac{d}{n}
\right\} \bigcup \left\{ S_{n,d}\right\} 
$$
 is equivariantly perfect. 
This means that we can modify our completeness criteria, so 
that for each pair of positive integers $(\n,\D)$ with $0<\n<n$ and
$\frac{\D}{\n} > \frac{d}{n}$ we require a set of relations
whose restriction in $\HG(S_{n_1,d_1})$ is zero when $d_1/n_1<\D/\n$ or
$d_1/n_1=\D/\n$ and $n_1<\n$, and when $(n_1,d_1) = (\n,\D)$ equals
the ideal of $\HG(S_{\n,\D})$ generated by the equivariant Euler class
$e_{\n,\D}$ of the normal to $S_{\n,\D}$ in $\C$. 

It is easy enough to prove that if $\gamma \in H_*^{\GG(\n,\D)}(\C(\n,\D)^{s})$ where $\D/\n > d/n$ and if
$r > n \n (g-1) + \n d - n \D$ then the image of the slant
product 
$
c_{r}(-\pi_!(\hat{\mathcal{V}}^* \otimes \mathcal{V})) \backslash \gamma
$
under the restriction map
$$
\HG(\C) \to \HG(S_{n_1,d_1}) \cong \HS_{\G(n_1,d_1)}(\C(n_1,d_1)^{ss})
\otimes \HS_{\G(n-n_1,d-d_1)}(U(n_1,d_1))
$$
is zero when $d_1/n_1 < \D/\n,$ and also when $d_1/n_1 = \D/\n$ and
$n_1<\n$.

\indent By Lefschetz duality, since $\C(\n,\D)^s
/\overline{\G}(\n,\D)= \mathcal{M}^s (\n,\D)$ is a manifold of dimension 
$$
D(\n,\D) = 2[(\n^2-1)(g-1) + g] = 2(\n^2 g - \n^2 +1)
$$
we have a natural map
$$
LD:H_*^{\GG(\n,\D)} (\C(\n,\D)^s) \cong H_*(\mathcal{M}^s(\n,\D)) \to 
H_{\GG(\n,\D)}^{D(\n,\D)-*}(\C(\n,\D)^{ss})
$$
such that if $\gamma \in H^{\GG(\n,\D)}_* (\C(\n,\D)^s)$ then $LD(\gamma)$ 
lies in the dual of $H^{\GG(\n,\D)}_{D(\n,\D)-*}(\C(\n,\D)^{ss})$ 
and takes a $\GG(\n,\D)$-equivariant  cycle on $\C(\n,\D)^{ss}$ to its intersection, 
modulo $\GG(\n,\D)$, with $\gamma$. 
When $\n$ and $\D$ are coprime then $\C(\n,\D)^{ss}$ equals
$\C(\n,\D)^s$ and its quotient by $\GG(\n,\D)$, namely $\M$, is a compact
manifold. In this case Lefschetz duality is essentially
Poincar\'{e} duality and the map $LD$ is an isomorphism. 

We need to consider the restriction of a relation of the
form $c_{r}(-\pi_!(\hat{\mathcal{V}}^* \otimes \mathcal{V})) \backslash \gamma$ to
$\HG(S_{\n,\D})$.
If $\gamma \in H_*^{\GG(\n,\D)}(\C(\n,\D)^s)$ where $\D/\n > d/n$ and
if $r = n\n(g-1)+ \n d - n \D +1
+j$,
then it turns out that the image of the slant product 
 $c_{r}(-\pi_!(\hat{\mathcal{V}}^* \otimes \mathcal{V})) \backslash \gamma$ under
the restriction map 
$$
\HG(\C) \to \HG(S_{\n,\D}) \cong \HS_{\G(\n,\D)}(\C(\n,\D)^{ss})
\otimes \HS_{\G(n-\n,d-\D)}(U(\n,\D))
$$
is the product 
$
(-a_1^{(1)})^j LD(\gamma) e_{\n,\D}
$
of the equivariant Euler class $e_{\n,\D}$ of the normal bundle to
$S_{\n,\D}$ in $\C$ with the image of $\gamma$ under the Lefschetz
duality map
$
LD$
and $j$ copies of minus the generator $a_1^{(1)} \in H^*_{\G(\n,\D)}(\C(\n,\D)^{ss})
  \cong \HS(BU(1)) \otimes \HS_{\GG(\n,\D)}(\C(\n,\D)^s)$ which comes from the
copy of the polynomial ring $\HS(BU(1))$. 
The proof of this is based on 
Porteous's Formula (as in Beauville's alternative proof \cite{B} 
of the theorem of Atiyah and Bott that the classes $\ar, \bjr, \fr$ 
generate $\HS(\mnd)$; cf. \cite{T2} and \cite{ES}), which allows us
to deduce that the Poincar\'{e} dual of $\Delta^s \times U(\n,\D)$ in
$$
\HS_{\G(\n,\D)}(\C(\n,\D)^s) \otimes  \HS_{\G(\n,\D)}(\C(\n,\D)^{ss}) \otimes 
\HS_{\G(n-\n,d-\D)}(U(\n,\D))
$$
is $ c_{\n^2 (g-1) +1}(-\pi_!(\hat{\mathcal{V}}^* \otimes \mathcal{V}_1).
$
In other words the restriction of $c_{\n^2(g-1)+1}(-\pi_!(\hat{\mathcal{V}}^* \otimes \mathcal{V}_1))$ 
to $\C(\n,\D)^s \times \C(\n,\D)^{ss} \times U(\n,\D)$ is the image
of $1$ under the $\G(\n,\D) \times \G(\n,\D) \times
\G(n-\n,d-\D)$-equivariant Thom-Gysin map 
associated to the inclusion of $\Delta^s \times U(\n,\D))$ in
$\C(\n,\D)^s \times
\C(\n,\D)^{ss} \times U(\n,\D))$. 
We can express the higher Chern classes of $-\pi_!(\hat{\mathcal{V}}^*
 \otimes \mathcal{V}_1)$ in a similar way \cite{EK} by using Fulton's Excess Porteous
 Formula \cite{F}.

Recall that given $\eta \in H^*_{\G(\n,\D)}(\C(\n,\D)^{ss}) \otimes
H^*_{\G(n-\n,d-\D)}(\C(n-\n,d-\D))$ we wish to find a relation whose
restriction to $\C(n_1,d_1)^{ss} \times U(n_1,d_1)$ is zero when
$d_1/n_1 < \D/\n$ or $d_1/n_1 = \D/\n$ and $n_1< n$, and when
$(n_1,d_1) = (\n,\D)$ equals $\eta e_{\n,\D}$. We have found such
a relation when $\eta$ lies in the image
of the Lefschetz duality map $LD$ which maps
$H_*^{\GG(\n,\D)}(\C(\n,\D)^s)$ to $H^*_{\GG(\n,\D)}(\C(\n,\D)^{ss})$ and
thus into 
$$
H^*_{\G(\n,\D)}(\C(\n,\D)^{ss}) = H^*_{\GG(\n,\D)}(\C(\n,\D)^{ss}) \otimes
\HS(BU(1)),
$$
and more generally when $\eta$ has the form
$
\eta = (-a_1^{(1)})^j LD(\gamma),
$
for some element $\gamma$ of $ H_*^{\GG(\n,\D)}(\C(\n,\D)^s)$. When $\n$ and $\D$ are coprime this
gives us all $\eta \in \HS_{\G(\n,\D)}(\Css)$. Moreover
$$
\HS_{\G(\n,\D)}(\C(\n,\D)^{ss}) \otimes \HS_{\G(n-\n,d-\D)}(\C(n-\n,d-\D))
$$
is generated as a module over $\HG(\C)$ by
$\HS_{\G(\n,\D)} (\C(\n,\D)^{ss})$,  so when $\n$ and $\D$ are coprime, we have now obtained the relations we
need from  the slant products
$$
\{ c_r(-\pi_!(\hat{\mathcal{V}}^* \otimes \mathcal{V})) \backslash \gamma: 
r \geq n \n (g-1) + \n d - n \D + 1, \gamma \in
H_*^{\GG(\n,\D)}(\C(\n,\D)^s) \},
$$
and a little more work reduces the range of $r$ and $\hat{d}$ (see \cite{EK} for more details).

This deals with the case when $\n$ and $\D$ are coprime, but the completeness criteria
have not yet been shown to hold for pairs $\n$ and $\D$ with common factors.
This was the original motivation for considering further modifications to the Yang--Mills
stratification. The difficulty with using the Yang--Mills stratification itself, or the stratification
of $\C$ by coarse type, is that in each case, although $n$ and $d$ are chosen to be coprime
so that semistability and stability coincide for $n$ and $d$, in the
construction of the stratification other $\n$ and $\D$ appear for which semistability and stability do not coincide.

\section{Stratifying a set of semistable points}
\renorm

In this section we shall describe briefly how to stratify the set $X^{ss}$ of semistable
points of a complex projective variety $X$ equipped with a linear action of a
complex reductive group $G$, so that the set $X^s$ of stable points of $X$ is
an open stratum (see \cite{G2,K2,MFK,N2} for background and \cite{K7} 
for more details).

We assume that $X$ has some stable points but also has semistable
points which are not stable. 
In \cite{K4,K6} it is described how one can blow $X$ up along a sequence of 
nonsingular  $G$-invariant subvarieties to obtain a $G$-invariant morphism 
$\tilde{X} \to X$ where $\tilde{X}$ is a complex projective variety acted on 
linearly by $G$ such that $\tilde{X}^{ss} = \tilde{X}^s$. 
The set $\tilde{X}^{ss}$ can be obtained from $X^{ss}$ as follows.  Let 
$r>0$ be the maximal dimension of a reductive subgroup of $G$ 
fixing a point of $X^{ss}$, and let $\mathcal{R}(r)$ be a set of representatives of conjugacy 
classes of all connected reductive subgroups $R$ of 
dimension $r$ in $G$ such that 
\[
Z^{ss}_{R} = \{ x \in X^{ss} :  \mbox{$R$ fixes $x$}\}
\]
is non-empty. Then
$
\bigcup_{R \in \mathcal{R}(r)} GZ^{ss}_{R}
$
is a disjoint union of nonsingular closed subvarieties of $X^{ss}$. The action of 
$G$ on $X^{ss}$ lifts to an action on the blow-up of 
$X^{ss}$ along $\bigcup_{R \in \mathcal{R}(r)} GZ_R^{ss}$ which can be linearised so that the complement 
of the set of semistable points in the blow-up is the proper transform of the 
subset $\phi^{-1}(\phi(GZ_R^{ss}))$ of $X^{ss}$ where $\phi:X^{ss} \to X/\!/G$ is the quotient 
map (see \cite{K4} 7.17). Moreover no semistable point in the blow-up is fixed by a 
reductive subgroup of $G$ of dimension at least $r$, and a semistable point in the blow-up 
is fixed by a reductive subgroup $R$ of 
dimension less than $r$ in $G$ if and only if it belongs to the proper transform of the 
subvariety $Z_R^{ss}$ of $X^{ss}$.

 If we repeat this process enough times, we obtain $\pi:\tilde{X}^{ss} \to X^{ss}$
 such that $\tilde{X}^{ss} = \tilde{X}^s$. Equivalently we can construct a sequence  
$$X_{(R_0)}^{ss} = X^{ss}, X_{(R_1)}^{ss},\ldots,X_{(R_\tau)}^{ss} = \tilde{X}^{ss}$$ 
where $R_1,\ldots,R_\tau$ are connected reductive subgroups of $G$ with 
$$r= \dim R_1 \geq \dim R_2 \geq \cdots \dim R_\tau \geq 1,$$ 
and if $1 \leq l \leq \tau$ then
 $X_{(R_l)}$ is the blow up of $X_{(R_{l-1})}^{ss}$ 
along its closed nonsingular subvariety $GZ_{R_l}^{ss}\cong G \times_{N_l} Z_{R_l}^{ss}$,
where $N_l$ is the normaliser of $R_l$ in $G$. Similarly, 
$\tilde{X}/\!/G = \tilde{X}^{ss}/G$ can be obtained from $X/\!/G$ by blowing 
up along the proper transforms of the images $Z_R /\!/N$ 
in $X/\!/G$ of the subvarieties $GZ_R^{ss}$ of $X^{ss}$ in decreasing order of $\dim R$.

If $1 \leq l \leq \tau$ then we have a $G$-equivariant stratification 
$$\{ \mathcal{S}_{\beta,l}: (\beta,l) \in \mathcal{B}_l\times \{ l \} \}$$ 
of $X_{(R_l)}$
by nonsingular 
$G$-invariant locally closed subvarieties such that one of the strata, indexed by $(0,l) \in \mathcal{B}_l \times\{ l \}$, 
coincides with the open 
subset $X_{(R_l)}^{ss}$ of $X_{(R_l)}$.  Here $\mathcal{B}_l$ is a finite subset of a fixed positive
Weyl chamber $\liet_+$ in the Lie algebra $\liet$ of a maximal compact torus $T$ of $G$. In fact
$\beta \in \liet_+$ lies in $\mathcal{B}_l$ if and only if $\beta$ is the closest point to 0 in the convex
hull in $\liet$ of some nonempty subset of the set of weights $\{ \alpha_0, \ldots, \alpha_n \}$ for the
linear action of $T$ on $X_{(R_l)}$.

There is a partial ordering on $\mathcal{B}_l$ given by $\gamma > \beta$ if 
$|\! | \gamma |\! | > |\! | \beta |\! |$, with $0$ as its 
minimal element, such that if 
$\beta \in \mathcal{B}_l$ then the closure in $X_l$ of the stratum $\mathcal{S}_{\beta,l}$ satisfies
\begin{equation} \label{sb}
\overline{\mathcal{S}_{\beta,l}} \subseteq \bigcup_{\gamma \in \mathcal{B}_l, \gamma \geq \beta} \mathcal{S}_{\gamma,l}.
\end{equation}
If $\beta \in \mathcal{B}_l$ and $\beta \neq 0$ then the stratum $\mathcal{S}_{\beta,l}$ retracts $G$-equivariantly 
onto its (tranverse) 
intersection with the exceptional divisor $E_l$ for the blow-up $X_{(R_l)} \to X_{(R_{l-1})}^{ss}$. This 
exceptional divisor is isomorphic to the projective bundle $\PP(\mathcal{N}_l)$ over $G\hat{Z}_{R_l}^{ss}$,
where $\hat{Z}_{R_l}^{ss}$ is the proper transform of $Z_{R_l}^{ss}$ in $X_{(R_{l-1})}^{ss}$ and
$\mathcal{N}_l$ is 
the normal bundle to $G \hat{Z}_{R_l}^{ss}$ in $X_{R_{l-1}}^{ss}$. The stratification 
$\{\mathcal{S}_{\beta,l}: \beta \in \mathcal{B}_l\}$ is determined 
by the action of $R_l$ on the fibres of $\mathcal{N}_l$ over $Z_{R_l}^{ss}$ (see \cite{K4} \S 7).


There is thus a stratification $\{ \Sigma_\gamma : \gamma \in \Gamma \}$ of $X^{ss}$ indexed by
\begin{equation}
\label{new*}
{\Gamma} = \{R_1\} \sqcup (\mathcal{B}_1 \backslash \{0\})\times\{1\} \sqcup  
 \cdots  
\sqcup  \{R_\tau\} \sqcup (\mathcal{B}_\tau \backslash \{0\})\times \{ \tau \} \sqcup \{ 0 \}
\end{equation}
defined as follows. We take as the highest stratum $\Sigma_{R_1}$ 
the nonsingular closed subvariety $GZ_{R_1}^{ss}$ whose 
complement in $X^{ss}$ can be naturally identified with 
the complement $X_{(R_1)}\backslash E_1$ of the exceptional divisor
$E_1$ in $X_{(R_1)}$.
We have $GZ^{ss}_{R_1} \cong G \times_{N_1} Z^{ss}_{R_1}$ where $N_1$ is the
normaliser of $R_1$ in $G$, and $Z^{ss}_{R_1}$ is equal to the set of semistable
points for the action of $N_1$, or equivalently for the induced action of $N_1/R_1$,
on $Z_{R_1}$, which is a union of connected components of the fixed point set of
$R_1$ in $X$. Moreover since $R_1$ has maximal dimension among those
reductive subgroups of $G$ with fixed points in $X^{ss}$, we have
$Z_{R_1}^{ss} = Z_{R_1}^s$
where $Z_{R_l}^s$ denotes the set of stable points for the action of $N_l/R_l$
on $Z_{R_l}$ for $1 \leq l \leq \tau$.

Next we take as strata the nonsingular locally closed subvarieties
$$ \Sigma_{\beta,1} = \mathcal{S}_{\beta,1}\backslash E_1 \mbox{ for } \beta \in \mathcal{B}_1 
\mbox{ with } \beta \neq 0
$$
of $X_{(R_1)}\backslash E_1 = X^{ss}\backslash GZ_{R_1}^{ss}$, whose complement 
in $X_{(R_1)}\backslash E_1$ is just 
$X_{(R_1)}^{ss}\backslash E_1=X_{(R_1)}^{ss}\backslash E_1^{ss}$ where $E_1^{ss} = X_{(R_1)}^{ss} \cap E_1$,
and then we take the intersection of 
$X_{(R_1)}^{ss}\backslash E_1$ with $GZ_{R_2}^{ss}$. This intersection is $G Z_{R_2}^s$ where
$Z_{R_2}^s$ is the set of stable points for the action of $N_2/R_2$ on $Z_{R_2}$, and its
complement in $X_{(R_1)}^{ss}\backslash E_1$ can be naturally identified with the 
complement in $X_{(R_2)}$ of the union of $E_2$ and the 
proper transform $\hat{E}_1$ of $E_1$.

The next strata are the nonsingular locally closed subvarieties
$$ \Sigma_{\beta,2} = \mathcal{S}_{\beta,2} \backslash (E_2 \cup \hat{E}_1) \mbox{ for } \beta \in \mathcal{B}_2
\mbox{ with } \beta \neq 0$$
of $X_{(R_2)}\backslash(E_2 \cup \hat{E}_1)$, whose complement in $X_{(R_2)}\backslash(E_2 \cup \hat{E}_1)$ is 
$X_{(R_2)}^{ss}\backslash(E_2 \cup \hat{E}_1)$, and the stratum after these 
is $GZ_{R_3}^s$. 
We repeat this process for $1 \leq l \leq \tau$ and take $X^s$ as our final stratum indexed by 0.

The given partial 
orderings on $\mathcal{B}_1,\ldots,\mathcal{B}_\tau$ together with the ordering in the expression (\ref{new*}) above
for $\Gamma$
induce a partial ordering on $\Gamma$, with $R_1$ as the maximal element
and 0 as the minimal element, such that the closure in $X^{ss}$ of the 
stratum $\Sigma_\gamma$ indexed by $\gamma \in \Gamma$ satisfies
\begin{equation} \label{stratclos}
\overline{\Sigma_\gamma} \subseteq \bigcup_{\tilde{\gamma}\geq \gamma } \Sigma_{\tilde{\gamma}}.
\end{equation}
It is possible to describe the strata $\Sigma_{\gamma}$ in more detail.
 Either $\Sigma_\gamma$ is
$GZ^s_{R_l},$ for some $l$, 
or else it is
$$\mathcal{S}_{\beta,l} \backslash (E_l \cup \hat{E}_{l-1} \cup ... \cup \hat{E}_1)$$
for some $l$ and $\beta \in \mathcal{B}_l \backslash \{ 0 \}$, in which case by \cite{K2} $\S$5 we have 
\begin{equation} \label{2.2} \mathcal{S}_{\beta,l} = G Y^{ss}_{\beta,l}
 \cong G \times_{P_{\beta}} Y^{ss}_{\beta,l} \end{equation}
where $Y^{ss}_{\beta,l}$ fibres over $Z^{ss}_{\beta,l}$  with fibre $\CC^{m_{\beta,l}}$
for some $m_{\beta,l} >0$, and $P_{\beta}$ is a parabolic subgroup of
$G$ with the stabiliser $\stab(\beta)$ of $\beta$ under the adjoint action of $G$
as its maximal reductive subgroup.  Here the fibration
$p_{\beta}:Y^{ss}_{\beta,l} \to Z^{ss}_{\beta,l}$ sends $y$ to a limit point of its orbit
under the complex one-parameter subgroup of $R_l$ generated by $\beta$. Moreover
\begin{equation} \label{2.3} \mathcal{S}_{\beta,l}\cap E_l = G (Y^{ss}_{\beta,l}\cap E_l) \cong G 
\times_{P_{\beta}} (Y^{ss}_{\beta,l} \cap E_l) \end{equation}
where $Y^{ss}_{\beta,l}\cap E_l$ fibres over $Z^{ss}_{\beta,l}$ with fibre
$\CC^{m_{\beta,l}-1}$ (see \cite{K4} Lemmas 7.6 and 7.11). Thus
\begin{equation} \label{2.4} \mathcal{S}_{\beta,l} \backslash E_l  \cong 
G \times_{P_{\beta}} (Y^{ss}_{\beta,l} \backslash E_l)
\end{equation}
where $Y^{ss}_{\beta,l} \backslash E_l$ fibres over $Z^{ss}_{\beta,l}$ with fibre
$\CC^{m_{\beta,l}-1} \times (\CC \backslash \{ 0 \})$. Furthermore if 
$\pi_l: E_l \cong \PP(\mathcal{N}_l) \to G \hat{Z}^{ss}_{R_l}$
is the projection then Lemma 7.9 of \cite{K4} tells us that when $x \in \hat{Z}^{ss}_{R_l}$ 
the intersection of $\mathcal{S}_{\beta,l}$ with the fibre $\pi_l^{-1}(x) = \PP(\mathcal{N}_{l,x})$
of $\pi_l$ at $x$ is the union of  those strata indexed by points in the adjoint orbit
$Ad(G)\beta$ in the stratification of $\PP(\mathcal{N}_{l,x})$ induced by the representation
$\rho_l$ of $R_l$ on the normal $\mathcal{N}_{l,x}$ to $G\hat{Z}^{ss}_{R_l}$ at $x$. 
Careful analysis \cite{K7} shows that we can, if we wish, 
replace the indexing set $\mathcal{B}_l \backslash \{ 0 \}$, whose elements
correspond to the $G$-adjoint orbits $Ad(G)\beta$ of elements of the indexing set for the
stratification of $\PP(\mathcal{N}_{l,x})$ induced by the representation $\rho_l$, by the set of their
$N_l$-adjoint orbits $Ad(N_l)\beta$. Then we still have (\ref{2.2}) -- (\ref{2.4}), but now if $q_\beta : P_\beta \to \stab(\beta)$ is the projection we have 
\begin{equation} \label{sgqy}
\Sigma_\gamma = \Sigma_{\beta,l} = \mathcal{S}_{\beta,l} 
 \backslash (E_l \cup \hat{E}_{l-1} \cup ... \cup \hat{E}_1)
\cong G \times_{Q_{\beta,l}} Y^{\backslash E}_{\beta,l} \end{equation}
where $Q_{\beta,l} = q_\beta^{-1} (N_l \cap \stab(\beta))$ and 
$$Y^{\backslash E}_{\beta,l} = Y^{ss}_{\beta,l} \backslash (E_l \cup \hat{E}_{l-1} \cup ... \cup \hat{E}_1)
 \cap p_\beta^{-1} \left(
Z^{ss}_{\beta,l} \cap \pi^{-1}_l(\hat{Z}^{ss}_{R_l}) \right),$$
and 
$p_\beta: Y^{\backslash E}_{\beta,l} \to Z^{ss}_{\beta,l} \cap \pi^{-1}_l(\hat{Z}^{ss}_{R_l}) $
is a fibration with fibre $\CC^{m_{\beta,l}-1} \times (\CC\backslash \{0\})$.

This process gives us a stratification 
$ \{ \Sigma_{\gamma}: \gamma \in \Gamma \} $
of $X^{ss}$
such that 
the stratum indexed by the minimal element 0 of $\Gamma$ coincides with the open subset $X^s$ of
$X^{ss}$. We shall apply this construction to obtain a stratification of $\mathcal{C}^{ss}$, and thus
inductively to refine the Yang-Mills stratification $\{\mathcal{C}_\mu:\mu \in \mathcal{M}\}$ of $\C$ by Harder--Narasimhan
type.

\section{Direct sums of stable bundles of equal slope}

\renorm

 In the good case when $n$ and $d$ are coprime, then
$\mathcal{C}^{ss} = \mathcal{C}^s$ and $\mnd = \mathcal{C}^{ss} / \G_c$ is a nonsingular
projective variety.
When $n$ and $d$ are not coprime, we can use the description of $\mnd$ as the geometric invariant theoretic
quotient $\mathcal{C}/\!/ \mathcal{G}_c$ to construct a partial desingularisation
$\tilde{\mathcal{M}}(n,d)$ of $\mnd$. From this construction we can use the
methods described in $\S$2 to obtain a stratification of $\Css$ with $\C^s$ as an open 
stratum, and thus (by considering the subquotients of the Harder--Narasimhan filtration) obtain
a stratification 
 of $\mathcal{C}$ refining the stratification
$\{ \mathcal{C}_\mu: \mu \in \mathcal{M} \}$. To understand this refined stratification
we need to use the description in \cite{K5} of the partial desingularisation
$\tilde{\mathcal{M}}(n,d)$. 

In fact in \cite{K5} $\tilde{\mathcal{M}}(n,d)$ is not constructed using the representation 
of $\mnd$ as the geometric invariant theoretic quotient of $\mathcal{C}$ by 
$\mathcal{G}_c$, although (as is noted at \cite[p.246]{K5}) this representation of $\mnd$ 
would lead to the same partial desingularisation. Instead in \cite{K5} $\mnd$ is represented 
as a geometric invariant theoretic quotient of a finite-dimensional nonsingular quasi-projective 
variety $R(\n,\D)$ by a linear action of $SL(p;\CC)$ where $p=\D+\n(1-g)$ with $\D >> 0$. We 
may assume that $\D >> 0$, since tensoring by a line bundle of degree $l$ gives an isomorphism 
of $\mnd$ with $\mathcal{M}(\n,\D+ \n l)$ for any $l \in {\bf Z}$. By \cite[Lemma 5.2]{N2} if $E$ 
is a semistable bundle over $\Sigma$ of rank $\n$ and degree $\D > \n(2g-1)$ where $g$ is the 
genus of $\Sigma$, then $E$ is generated by its sections and $H^1(\Sigma,E)=0$. If $p = \D + \n(1-g)$, 
this implies that $\dim H^0(\Sigma,E)=p$ and that there is a holomorphic map $h$ from $\Sigma$ to 
the Grassmannian $G(\n,p)$ of $\n$-dimensional quotients of $\CC^p$ such that the pullback $E(h)=h^*T$ 
of the tautological bundle $T$ on $G(\n,p)$ is isomorphic to $E$.

Let $R(\n,\D)$ be the set of all holomorphic maps $h:\Sigma \to G(\n,p)$ such that $E(h)=h^*T$ has 
degree $d$ and the map on sections $\CC^p \to H^0(\Sigma,E(h))$ induced from the quotient bundle 
map $\CC^p \times \Sigma \to E(h)$ is an isomorphism. For $\D >>0$ this set $R(\n,\D)$ has
the structure of a nonsingular 
quasi-projective variety and there is a quotient $\mathcal{E}$ of the trivial bundle of rank $p$ over $R(\n,\D)\times \Sigma$ 
satisfying the following properties (see \cite[$\S 5$]{N2}).

\indent (i) If $h \in R(\n,\D)$ then the restriction of $E$ to $\{h\} \times \Sigma$ is the pullback $E(h)$ of the tautological 
bundle $T$ on $G(\n,p)$ along the map $h:\Sigma \to G(\n,p)$.\\
\indent (ii) If $h_1$ and $h_2$ lie in $R(\n,\D)$ then $E(h_1)$ and $E(h_2)$ are isomorphic as bundles over $\Sigma$ 
if and only if $h_1$ and $h_2$ lie in the same orbit of the natural action of $GL(p;\CC)$ on $R(\n,\D)$.\\
\indent (iii) If $h \in R(\n,\D)$ then the stabiliser of $h$ in $GL(p;\CC)$
is isomorphic to the group $\mbox{Aut}(E(h))$ 
of complex analytic automorphisms of $E(h)$.

If $N >> 0$ then $R(\n,\D)$ can be embedded as a quasi-projective subvariety of the product $(G(\n,p))^N$ by a map of the form
$
h \mapsto (h(x_1),.\ldots,h(x_N))
$
where $x_1,\ldots,x_N$ are points of $\Sigma$. This embedding gives a linearisation of the action of $GL(p;\CC)$ 
on $R(\n,\D)$. If $N >>0$ and $d >>0$ then we also have the following:\\
\indent (iv) The point $h \in R(\n,\D)$ is semistable in the sense of geometric invariant theory for the linear action of 
$SL(p;\CC)$ on $R(\n,\D)$ if and only if $E(h)$ is a semistable bundle. If $h_1$ and $h_2$ lie in $R(\n,\D)^{ss}$ 
then they represent the same point of $R(\n,\D)/\!/SL(p;\CC)$ if and only if $\mbox{gr}(E(h_1)) \cong \mbox{gr}(E(h_2))$,
and thus there is a natural identification of $\mnd$ with $R(\n,\D)/\!/SL(p;\CC)$ (see for example \cite[$\S$5]{N2}).

It is shown in \cite{K3} that the Atiyah--Bott formulas for the equivariant Betti numbers of 
$\mathcal{C}^{ss}$ can be obtained by stratifying $R(n,d)$ instead of $\C(n,d)$, and in 
fact throughout this paper we could work with either $R(n,d)$ or $\C(n,d)$. In particular properties (i) to (iv) above imply 
that the analysis in \cite{K5} of the construction of the partial desingularisation of $\tilde{\mathcal{M}}(\n,\D)$ as $\tilde{R}(\n,\D)/\!/SL(p;\CC)$ 
applies equally well if we work with $\tilde{\C}/\!/\mathcal{G}_c$.

To describe the construction of $\tilde{\mathcal{M}}(n,d)$, first of all we need to find a set $\mathcal{R}$ of representatives of the
conjugacy classes of reductive subgroups of $SL(p;\CC)$ which occur as the connected components of stabilisers of semistable
points of $R(n,d)$. In fact it is slightly simpler to describe the corresponding subgroups of $G=GL(p;\CC)$, and since the
central one-parameter subgroup of $GL(p;\CC)$ consisting of nonzero scalar multiples of the identity acts
trivially on $R(n,d)$, finding stabilisers in $GL(p;\CC)$ is essentially equivalent to finding stabilisers in $SL(p;\CC)$.
By \cite[pp. 248-9]{K5} such conjugacy classes in $GL(p;\CC)$ correspond to unordered sequences
$(m_1,n_1),...,(m_q,n_q)$ of pairs of positive integers such that
$m_1 n_1 + ... + m_q n_q = n$
and $n$ divides $n_i d$ for each $i$. A representative $R$ of the corresponding conjugacy class is given
by
$$R= GL(m_1;\CC) \times ... \times GL(m_q;\CC)$$
embedded in $GL(p;\CC)$ using a fixed isomorphism of $\CC^p$ with $\bigoplus_{i=1}^q (\CC^{m_i} \otimes \CC^{p_i}))$
where $d_i = n_i d/n$ and $p_i = d_i + n_i(1-g) = n_i p/n$.  Then $GZ_R^{ss}$ consists of all those holomorphic structures $E$ with
\begin{equation}
\label{newo*}
E \cong (\CC^{m_1} \otimes D_1) \oplus \cdots \oplus (\CC^{m_p} \otimes D_q)
\end{equation}
where $D_1, \ldots ,D_q$ are all semistable and $D_i$ has rank $n_i$ and degree $d_i$, while
$GZ_R^{s}$ consists of all those holomorphic structures $E$ with
\begin{equation}
E \cong (\CC^{m_1} \otimes D_1) \oplus \cdots \oplus (\CC^{m_p} \otimes D_q)
\end{equation}
where $D_1, \ldots ,D_q$ are all stable and not isomorphic to one another,
and $D_i$ has rank $n_i$ and rank $d_i$.  Moreover the normaliser $N$ of $R$ in
$GL(p; \CC)$ has connected component
\begin{equation} \label{N_0} N_0 \cong \prod_{1 \leq i \leq q} (GL(m_i;\CC) \times GL(p_i;\CC))/\CC^*
\end{equation}
where $\CC^*$ is the diagonal central one-parameter subgroup of $GL(m_i;\CC) \times GL(p_i;\CC)$,
and $\pi_0(N) = N/N_0$ is the product
\begin{equation} \label{N/N_0} \pi_0(N) = \prod_{j\geq 0, k \geq 0} Sym({\#\{i:m_i=j \mbox{ and }n_i=k\}})
\end{equation}
where $Sym(b)$ denotes the symmetric group of permutations of a set with $b$ elements.
Furthermore in terms of the notation of \S 2, if $R=R_l$ then a holomorphic structure
belongs to one of the strata $\Sigma_{\beta,l}$ with $\beta \in \mathcal{B}_l \backslash \{ 0 \}$
if and only if it has a filtration 
$0=E_0 \subset E_1 \subset ... \subset E_s= E$
such that $E$ is not isomorphic to $\bigoplus_{1 \leq k \leq s} E_k/E_{k-1} $ but
$$\bigoplus_{1 \leq k \leq s} E_k/E_{k-1} 
\cong (\CC^{m_1} \otimes D_1) \oplus \cdots \oplus (\CC^{m_p} \otimes D_q)
$$
where $D_1, \ldots ,D_q$ are all stable and not isomorphic to one another,
and $D_i$ has rank $n_i$ and rank $d_i$ \cite[p. 248]{K5}. Thus to understand the
refined Yang-Mills 
stratification   
  of $\mathcal{C}$, we need to
study refinements
$$0=E_0 \subset E_1 \subset ... \subset E_s= E$$
of the Harder-Narasimhan filtration of a holomorphic bundle $E$, such that
each subquotient $E_j/E_{j-1}$ is a direct sum of stable bundles all of the same slope.
For this recall the following standard result 
(cf. the proof of \cite[Lemma 3.2]{K5} and \cite{Se}).

\begin{prop}
\label{nineteen}
Any semistable bundle $E$ has a canonical subbundle of the form 
\[
(\CC^{m_1} \otimes D_1) \oplus \cdots \oplus (\CC^{m_q}
\otimes D_q)
\]
with $D_1,\ldots,D_q$ not isomorphic to each other and all stable of
the same slope as $E$, such that any other subbundle of the same form 
\[
(\CC^{m'_1} \otimes D'_1) \oplus \cdots \oplus (\CC^{m'_r}
\otimes D'_r) ,
\]
with $D'_1,...,D'_r$ not isomorphic to each other and all stable of
the same slope as $E$, satisfies  $r \leq q$ and (after
permuting the order of $D_1,...,D_q$ suitably) $D_j \cong
D'_j$ for all $1 \leq j \leq r$ and the inclusion
\[
(\CC^{m'_1} \otimes D_1) \oplus \cdots \oplus (\CC^{m'_r}
\otimes D_r) \to E 
\]
factors through the inclusion 
$ (\CC^{m_1} \otimes D_1) \oplus \cdots \oplus (\CC^{m_q}
\otimes D_q) \to E $
via injections $\CC^{m'_j} \to \CC^{m_j}$ for $1 \leq j \leq
r.$
\end{prop} 

If we choose a subbundle of $E$ of maximal rank among those of the
required form, this follows immediately from

\begin{lem} \label{e4.2} 
Let 
$E,D_1,\ldots,D_q,D'_1,\ldots,D'_r$ be bundles over $\Sigma$ all of
the same slope, with $E$ semistable, $D_1,...,D_q$ and $D'_1,\ldots,D'_r$  all stable
and $D_j = D'_j \mbox{ for } 1 \leq j \leq k$ for some $0 \leq k \leq \min\{q,r\}$,
but with no other isomorphisms between the bundles $D_1,...,D_q$ and
$D'_1,...,D'_r$. 
If 
$$\alpha: (\CC^{m_1} \otimes D_1) \oplus \cdots \oplus (\CC^{m_q}
\otimes D_q) \to E  $$
and 
\[
\beta : (\CC^{m'_1} \otimes D'_1) \oplus \cdots \oplus (\CC^{m'_r}
\otimes D'_r) \to E 
\]
are injective bundle homomorphisms, then there
exist nonnegative integers $n_1,...,n_k$ and linear injections 
$ i_j: \CC^{m_j} \to \CC^{n_j}$ and $i'_j: \CC^{m'_j} \to \CC^{n_j}$
for $1 \leq j \leq k$ and an injective bundle homomorphism
$\gamma$ from
$$ (\CC^{n_1} \otimes D_1) \oplus \cdots \oplus
(\CC^{n_k} \otimes D_k)  \oplus 
(\CC^{m_{k+1}} \otimes D_{k+1}) \oplus \cdots \oplus
(\CC^{m_q} \otimes D_q) $$
$$ \oplus 
(\CC^{m'_{k+1}} \otimes D'_{k+1}) \oplus \cdots \oplus
(\CC^{m'_r} \otimes D'_r) $$
to $E$ 
such that $\alpha$ and $\beta$ both factorise through $\gamma$ via the
injections $i_j$ and $i'_j$ for $1 \leq j \leq k$ in the obvious way.
\end{lem}
{\bf Proof:} Consider the kernel of 
\begin{eqnarray*}
\alpha \oplus \beta: (\CC^{m_1+m'_1} \otimes D_1) \oplus \cdots
\oplus (\CC^{m_k+m'_k} \otimes D_k) & \oplus & \\
(\CC^{m_{k+1}} \otimes D_{k+1}) \oplus \cdots \oplus
(\CC^{m_q} \otimes D_q) & \oplus &\\
(\CC^{m'_{k+1}} \otimes D'_{k+1}) \oplus \cdots \oplus
(\CC^{m'_r} \otimes D'_r) & \to & E.
\end{eqnarray*}
The proof of \cite[Lemma 5.1 (iii)]{N2} shows that this kernel 
is a subsheaf of the domain of $\alpha \oplus \beta$ which has the
same slope as $E, D_1,...,D_q,$ and $D'_1,...,D'_r$.
Induction on $m_1 + ... +m_q + m'_1 +...+m'_r$
using \cite[Lemma 5.1]{N2} and the obvious projection from
$\CC^{m_1}$ onto $\CC^{m_1 -1}$ shows that
such a subsheaf is in fact of the form 
\[
(U_1 \otimes D_1) \oplus \cdots \oplus (U_q \otimes D_q) \oplus 
(U_{q+1} \otimes D'_{k+1}) \oplus \cdots \oplus (U_{q-k+r} \otimes D'_r)
\]
for some linear subspaces $U_j$ of $\CC^{m_1+m'_1},...,
\CC^{m_{k+1}},...,\CC^{m'_r}$. Then since $\alpha$ and
$\beta$ are injective, we must have
$
\mbox{ker}(\alpha \oplus \beta) = (U_1 \otimes D_1) \oplus \cdots
\oplus (U_k \otimes D_k)
$
for linear subspaces $U_j$ of $
\CC^{m_j} \oplus \CC^{m'_j}$ satisfying 
$U_j \cap \CC^{m_j} = \{0\} = 
U_j \cap \CC^{m'_j}$
for $1 \leq j \leq k$.  The result follows easily.

\bigskip


A similar argument gives us

\begin{cor} \label{lem34}
Let $E$ be semistable and $D'_1,\ldots, D'_r$ all stable of the same
slope as $E$ and not isomorphic to each other. Then 
\[
H^0(\Sigma,((\CC^{m'_1} \otimes D'_1) \oplus \cdots
\oplus (\CC^{m'_r} \otimes D'_r))^* \otimes E) \cong \bigoplus_{j=1}^k
((\CC^{m'_j})^* \otimes \CC^{m_j})
\]
where 
$
(\CC^{m_1} \otimes D_1) \oplus \cdots \oplus (\CC^{m_q} \otimes D_q)
$
is the canonical subbundle of this form associated to $E$ as in Proposition
\ref{nineteen}, and without loss of generality we assume $D_j
\cong D'_j$ for $1 \leq j \leq k$ for some $0 \leq k \leq \min\{q,r\}$
and that there are no other isomorphisms between the bundles 
$D_1,...,D_q$ and $D'_1,...,D'_r$.
\end{cor}

\begin{df} \label{gr}
With the notation above we set 
\[
\mbox{gr}(E) = \bigoplus_{i=1}^q \bigoplus_{j=1}^r \left( \CC^{m_{ij}} \otimes D_i \right).
\]
for any semistable bundle $E$.
\end{df}

\begin{rem} \label{jh}
Of course, by  the Jordan--H\"{o}lder theorem, given {\em any} filtration
$0=D_0 \subset D_1 \subset ... \subset D_t = E$
of $E$ such that $D_j/D_{j-1}$ is a direct sum of stable bundles of the same slope as
$E$, we have
$\bigoplus_{j=1}^{t}   D_j/D_{j-1} 
\cong \mbox{gr}(E).$
\end{rem}

\section{Maximal and minimal Jordan--H\"{o}lder filtrations}
\renorm

Recall that the Harder-Narasimhan filtration of a holomorphic bundle $E$ over $\Sigma$
is a canonical filtration
$$ 0 = F_0 \subset F_1 \subset ... \subset F_s = E$$
of $E$ such that $F_j/F_{j-1}$ is semistable and 
$\mbox{slope}(F_j/F_{j-1}) > \mbox{slope}(F_{j+1}/F_j)$ for $0 < j < s$. In 
the last section we saw
that any semistable bundle $E$ has a canonical maximal subbundle of the
form 
\[
(\CC^{m_1} \otimes D_1) \oplus \cdots \oplus (\CC^{m_q}
\otimes D_q) 
\]
where $D_1,\ldots,D_q$ are not isomorphic to each other and are all stable of
the same slope as $E$. This subbundle is nonzero if $E \neq 0$,
since any nonzero semistable bundle is either stable 
itself or it has a proper stable subbundle of the same slope.

Therefore 
any semistable bundle $E$ has a canonical filtration 
\begin{equation}
\label{new1}
0 = E_0 \subset E_1 \subset E_2 \subset \cdots \subset E_r = E
\end{equation}
whose subquotients are direct sums of stable bundles, which  is defined inductively so that 
\begin{equation}
\label{newd}
E_j/E_{j-1} \cong (\CC^{m_{1j}} \otimes D_1) \oplus \cdots \oplus (\CC^{m_{qj}} \otimes D_q)
\end{equation}
where $D_1,...,D_q$ are stable nonisomorphic bundles all of the same slope as $E$
with nonnegative integers  $m_{ij}$ for $1 \leq i \leq q$ and $1 \leq j \leq r$, 
and $E_j/E_{j-1}$ is the maximal subbundle of $E/E_{j-1}$ of this form. If, moreover, we assume that 
\[
\sum_{j=1}^{r} m_{ij} > 0 \mbox{ for all } 1 \leq i \leq q
\]
then the filtration (\ref{new1}), the bundles $D_i$ and integers $m_{ij}$ (for $1 \leq i \leq q$ and $1 \leq j \leq r$) 
and the decompositions (\ref{newd}) are canonically associated to $E$ up to isomorphism of the bundles $D_i$, the usual 
action of $GL(m_{ij};\CC)$ on $\CC^{m_{ij}}$ and the obvious action of the permutation group $Sym(q)$ on this data.
We can generalise the definition to the case when  
$E$ is not necessarily semistable, by applying this construction 
to the subquotients of the Harder-Narasimhan
filtration 
of $E$. This gives us a canonical refinement
$$0=E_0 \subset E_1 \subset ... \subset E_t = E$$
of the Harder-Narasimhan filtration such that each subquotient $E_j/E_{j-1}$ is the maximal
subbundle of $E/E_{j-1}$ which is a direct sum of stable bundles all having maximal slope 
among the nonzero subbundles of $E/E_{j-1}$. We shall call this
refinement of the Harder-Narasimhan filtration the maximal Jordan--H\"{o}lder 
filtration of $E$. 

\begin{df} \label{def81}
Let $s$, $q_1,...,q_s$ and $r_1,...,r_s$ be positive integers and let $d_k$, $n_{ik}$ and $m_{ijk}$ 
(for $1 \leq i \leq q_k$, $1 \leq j \leq r_k$
and $1 \leq k \leq s$) be integers satisfying
$$ n_{ik} > 0,  \quad m_{ijk} \geq 0, \quad  n_k  > 0$$
and 
$$ n = \sum_{k=1}^s n_k, 
\quad d = \sum_{k=1}^s d_k, \quad \frac{d_k n_{ik}}{n_k} = d_{ik} \in {\ZZ}, 
\quad \frac{d_1}{n_1} > \frac{d_2}{n_2} > ... > \frac{d_s}{n_s}$$
\end{df}
where
$$n_k = \sum_{i=1}^{q_k} \sum_{j=1}^{r_k} n_{ik}m_{ijk}.$$
Denote by $[{\bf d,n,m}] = [(d_k)_{k=1}^s,(n_{ik})_{i=1,k=1}^{q_k,s},
(m_{ijk})_{i=1,j=1,k=1}^{q_k,r_k,s}]$ the orbit of 
$$({\bf d,n,m}) = ((d_k)_{k=1}^s,(n_{ik})_{i=1,k=1}^{q_k,s},
(m_{ijk})_{i=1,j=1,k=1}^{q_k,r_k,s})$$ 
under the action of the product of symmetric
groups $\Sigma_{q_1} \times ... \times \Sigma_{q_k}$ on 
the set of such sequences, and let $\mathcal{I} = \mathcal{I}(\n,\D)$ denote the set of all such orbits, 
for fixed $\n$ and $\D$. Given
$ [{\bf d},{\bf n},{\bf m}] 
\in \mathcal{I}(\n,\D) $
let $s([{\bf d},{\bf n},{\bf m}] = s$ and
let $S^{maxJH}_{[{\bf d},{\bf n},{\bf m}]}$ denote the subset of $\mathcal{C}$ consisting of those holomorphic 
structures on our fixed smooth bundle of rank $\n$ and degree $\D$ whose maximal Jordan--H\"{o}lder
filtration 
\[
0 = E_{0,1} \subset E_{1,1} \subset \cdots \subset E_{r_1,1}=E_{0,2}\subset E_{1,2} \subset \cdots \]
\[ \cdots
\subset E_{r_{s-1}, s-1} = E_{r_s, 0} \subset \cdots \subset E_{r_s,s} = E
\]
satisfies
\[
E_{j,k}/E_{j-1,k} \cong (\CC^{m_{1jk}} \otimes D_{1k}) \oplus \cdots \oplus (\CC^{m_{q_k jk}} \otimes D_{q_k k})
\] 
for $1 \leq k \leq s$ and $1 \leq j \leq q_k$, where $D_{1 k},\ldots,D_{q_k k}$ are nonisomorphic stable bundles with 
\[
\mbox{rank}(D_{ik}) = n_{ik} \mbox{ and } \mbox{deg}(D_{ik}) = d_{ik}
\]
and $E_{j,k}/E_{j-1,k}$ is the maximal subbundle of $E/E_{j-1,k}$ isomorphic to a direct sum of stable 
bundles of  slope $d_k/n_k$. To make the notation easier on the eye, 
$S^{maxJH}_{[{\bf d},{\bf n},{\bf m}]}$ will often be denoted by $S_{[{\bf d},{\bf n},{\bf m}]}$.

 Let $\mathcal{I}^{ss}$ denote the subset of $\mathcal{I}$ consisting of all
orbits $[{\bf d},{\bf n},{\bf m}] $ for which $s([{\bf d},{\bf n},{\bf m}] = 1$. For simplicity we shall write
$[{\bf n},{\bf m}] $ for $[{\bf d},{\bf n},{\bf m}] $  when $s([{\bf d},{\bf n},{\bf m}])=1 $ (which means that
${\bf d} = (d)$).

We have now proved

\begin{lem}
\label{disj}
$\mathcal{C}$ is the disjoint union of the subsets
\[
\{ S_{[{\bf d},{\bf n},{\bf m}]  } : [{\bf d}, {\bf n},{\bf m}] \in \mathcal{I}\}.
\]
\end{lem}

\begin{rem}
\label{fdisj}
Note that when $s( [{\bf d},{\bf n},{\bf m}] )  =1$ and $q_1=r_1=1$ and $ n_{11} = \n$ and $m_{111}=1$ 
then we get $S_{[(n),(1)]}=\mathcal{C}^s$, and moreover
$$\C(n,d)^{ss} = \bigcup_{[{\bf d},{\bf n},{\bf m}] \in \mathcal{I}, s([{\bf d},{\bf n},{\bf m}])=1} S_{[{\bf d},{\bf n},{\bf m}] }
= \bigcup_{[{\bf n},{\bf m}] \in \mathcal{I}^{ss}} S_{[{\bf n},{\bf m}]} .$$
\end{rem}

\begin{rem} \label{4.3a}
Let $E$ be a semistable holomorphic structure on $\mathcal{E}$. If $E$ represents an element of the
closure of $S_{[{\bf n, m}]}$ in $\C(n,d)^{ss}$ for some 
$$ [{\bf n,m}] = [(n_i)_{i=1}^q, (m_{ij})_{i=1,j=1}^{q,r}]   \in \mathcal{I}^{ss}, $$ then $E$
has a filtration
$0 = E_0 \subset E_1 \subset E_2 \subset \cdots \subset E_r = E
$ such that if $1 \leq j \leq r$ then
$$
E_j/E_{j-1} \cong (\CC^{m_{1j}} \otimes D_1) \oplus \cdots \oplus (\CC^{m_{qj}} \otimes D_q)$$
where $D_1,...,D_q$ are semistable bundles all having the same slope as $E$, but
this filtration is not necessarily the maximal Jordan--H\"{o}lder filtration of $E$. If the
bundles $D_1,...,D_q$ are not all stable or if two of them are isomorphic to each other, then
$$ \mbox{dim Aut}(gr(E)) > \sum_{i=1}^q (\sum_{j=1}^r m_{ij})^2$$
and so $E$ lies in $S_{[\bf{n',m'}]}$ where ${\bf m'} = (m'_{ij})_{i=1,j=1}^{q',r'}$ satisfies
$$ \sum_{i=1}^{q'} (\sum_{j=1}^{r'} m'_{ij})^2 >  \sum_{i=1}^q (\sum_{j=1}^r m_{ij})^2.$$
If, on the other hand, $D_1,...,D_q$ are all stable and not isomorphic to each other, then
the maximal Jordan--H\"{o}lder filtration of $E$ is of the form
$$0 = E'_0 \subset E'_1 \subset E'_2 \subset \cdots \subset E'_r = E
$$ with
$$
E'_j/E'_{j-1} \cong (\CC^{m'_{1j}} \otimes D_1) \oplus \cdots \oplus (\CC^{m'_{qj}} \otimes D_q)$$
for $1 \leq j \leq r'$, where $1 \leq r' \leq r$ and 
$$m'_{i1} + ... + m'_{ir'} = m_{i1} + ...+ m_{ir}$$
and
$m'_{i1} + ... + m'_{ij} \geq  m_{i1} + ...+ m_{ij}$
for $1 \leq i \leq q$ and $1 \leq j \leq r'$. Thus we can define a partial order $\geq$ on $\mathcal{I}^{ss}$
such that $[{\bf n',m'}] \geq [{\bf n,m}]$ if and only if either 
$$ \sum_{i=1}^{q'} (\sum_{j=1}^{r'} m'_{ij})^2 > \sum_{i=1}^q (\sum_{j=1}^r m_{ij})^2$$
or ${\bf n}' = {\bf n}$ and $1 \leq r' \leq r$ and 
$$m'_{i1} + ... + m'_{ir'} = m_{i1} + ...+ m_{ir}$$
and
$m'_{i1} + ... + m'_{ij} \geq  m_{i1} + ...+ m_{ij}$
for $1 \leq i \leq q$ and $1 \leq j \leq r'$, and then the closure 
 of $S_{[{\bf n, m}]}$ in $\C(n,d)^{ss}$ is contained in
$$\bigcup_{[{\bf n',m'}] \geq [{\bf n,m}]} S_{[{\bf n', m'}]}.$$
Using (\ref{po}), we can then extend this partial order to $\mathcal{I}$  so that
$$\overline{ S_{[{\bf d,n,m}]}} \subseteq
\bigcup_{[{\bf d',n',m'}] \geq [{\bf d,n,m}]} S_{[{\bf d',n', m'}]}.$$
\end{rem}

\begin{prop}
\label{4.4}
Let $[{\bf n},{\bf m}] = [(n_i)^{q}_{i=1},(m_{ij})_{i=1,j=1}^{q,r}] \in \mathcal{I}^{ss}$.  
Then $S_{[{\bf n},{\bf m}]}$ is a locally closed complex submanifold of $\mathcal{C}^{ss}$ of finite 
codimension 
\[
\sum_{i=1}^{q}\sum_{j=1}^{r-1} m_{ij} m_{i j+1} + (g-1) \left( \sum_{i_1,i_2=1}^{q} \sum_{1 \leq j_1 \leq j_2\leq r} 
m_{i_1 j_1} m_{i_2 j_2} n_{i_1} n_{i_2} - \sum_{i=1}^{q} \sum_{j=1}^{r} (m_{ij}n_i)^2 \right) .
\]
\end{prop}

{\bf Proof:} (cf. \cite[$\S$7]{AB}) The rank and degree are the only $C^\infty$ invariants of 
a vector bundle over $\Sigma$. Thus we may choose a $C^\infty$ isomorphism of our fixed 
$C^\infty$ bundle $\mathcal{E}$ over $\Sigma$ with a bundle of the form 
\[
\bigoplus_{i=1}^{q} \bigoplus_{j=1}^{r} ({\bf C}^{m_{ij}} \otimes \mathcal{D}_i )
\]
where $\mathcal{D}_i$ is a fixed $C^\infty$ bundle over $\Sigma$ of rank $n_i$ and degree 
$d_i = n_i \D / \n$ for $1 \leq i \leq q$.

Let $\mathcal{Y}_{[{\bf n},{\bf m}]}$ be the subset of $\mathcal{C}^{ss}$ consisting of all semistable 
holomorphic structures $E$ on $\mathcal{E}$ for which the subbundles 
\[
E_j = \bigoplus_{i=1}^{q} \bigoplus_{k=1}^{j} ({\bf C}^{m_{ik}} \otimes \mathcal{D}_i)
\]
are holomorphic for $1 \leq j \leq r$, and for which there are nonisomorphic stable holomorphic structures 
$D_1,\ldots,D_q$ on the $C^\infty$ bundles $\mathcal{D}_1,\ldots,\mathcal{D}_q$ such that the 
natural identification of $E_j/E_{j-1}$ with 
$ \bigoplus_{i=1}^{q} {\bf C}^{m_{ij}} \otimes \mathcal{D}_i $
becomes an isomorphism of $E_j/E_{j-1}$ with 
$ \bigoplus_{i=1}^{q} {\bf C}^{m_{ij}} \otimes D_i $
for $1 \leq j \leq r$,  and finally for each $1 \leq j \leq r$ the quotient $E_j/E_{j-1}$ is 
the maximal subbundle of $E/E_{j-1}$ isomorphic to a direct sum of stable bundles 
of the same slope as $E/E_{j-1}$. Let $\G_c[{\bf n},{\bf m}]$ be the subgroup of the 
complexified gauge group $\mathcal{G}_c$ consisting of all $C^\infty$ complex automorphisms 
of $\mathcal{E}$ which preserve the filtration of $\mathcal{E}$ by the subbundles 
$\bigoplus_{i=1}^{q} \bigoplus_{k=1}^{j} ({\bf C}^{m_{ik}} \otimes \mathcal{D}_i)$
and the decomposition of $E_j/E_{j-1}$ as
$\bigoplus_{i=1}^{q} ({\bf C}^{m_{ij}} \otimes \mathcal{D}_i)$
up to the action of the general linear groups $GL(m_{ij};\CC)$ and the permutation
groups
$$Sym({\# \{i:n_i = k \mbox{ and } m_{ij} = l_j, j=1,...,r\} })$$
for all nonnegative integers $k$ and $l_1,...,l_r$.
Since the filtration (\ref{new1}) and decompositions of $E_1/E_0$,...,$E_r/E_{r-1}$ 
are canonically associated to $E$ up to the actions of these general linear 
and permutation groups, we have 
\[
S_{[{\bf n},{\bf m}]} = \mathcal{G}_c \mathcal{Y}_{[{\bf n},{\bf m}]} \cong  \mathcal{G}_c \times_{\G_c[{\bf n},{\bf m}]} 
\mathcal{Y}_{[{\bf n},{\bf m}]}.
\]
As in \cite[$\S$7]{AB} we have that $\mathcal{Y}_{[{\bf n},{\bf m}]}$ is an open subset of an affine subspace of the infinite-dimensional 
affine space $\mathcal{C}$ and the injection
\begin{equation} \label{inj}
\mathcal{G}_c \times_{\G_c[{\bf n},{\bf m}]} \mathcal{Y}_{[{\bf n},{\bf m}]} \to \mathcal{C}
\end{equation}
is holomorphic with image $S_{[{\bf n},{\bf m}]}$. 

If $E \in \mathcal{Y}_{[{\bf n},{\bf m}]}$, let $\mbox{End}'E$ be the subbundle of 
$\mbox{End}E$ consisting of holomorphic endomorphisms of $E$ preserving 
the maximal Jordan--H\"{o}lder filtration (\ref{new1}) and decomposition (\ref{newd}) 
up to isomorphism of the bundles $D_i$ and the vector spaces $\CC^{m_{ij}}$. 
Let $\mbox{End}''E$ be the quotient of $\mbox{End}E$ 
by $\mbox{End}'E$. The normal to the $\mathcal{G}_c$-orbit of $E$ in $\mathcal{C}$ can be 
canonically identified with $H^1(\Sigma,\mbox{End}E)$ (see \cite[$\S$7]{AB}) and 
the image of $T_E \mathcal{Y}_{[{\bf n},{\bf m}]}$ in this can be canonically identified 
with the image of the natural map
\[
H^1(\Sigma,\mbox{End}'E) \to H^1(\Sigma,\mbox{End}E)
\]
which fits into the long exact sequence of cohomology induced by the short exact sequence of bundles
\[
0 \to \mbox{End}'E \to \mbox{End}E \to \mbox{End}''E \to 0.
\]
Thus we get an isomorphism
\[
T_E \mathcal{C}/ (T_E \mathcal{Y}_{[{\bf n},{\bf m}]} + T_E \mathcal{O}) \cong H^1(\Sigma, \mbox{End}''E)
\]
where $\mathcal{O}$ is the $\mathcal{G}_c$-orbit of $E$ in $\mathcal{C}$.

We have short exact sequences
\[
0 \to (E/E_1)^* \otimes E \to \mbox{End}E \to E_1^* \otimes E \to 0
\]
and 
\[
0 \to E_1^* \otimes E_1 \to E_1^* \otimes E \to E_1^* \otimes (E/E_1) \to 0.
\]
Let $K$ be the kernel of the composition of surjections
\[
\mbox{End}E \to E_1^* \otimes E \to E_1^* \otimes (E/E_1).
\]
Then we have short exact sequences 
\[
0 \to K \to \mbox{End}E \to E_1^* \otimes (E/E_1) \to 0
\]
and
\[
0 \to (E/E_1)^* \otimes E \to K \to E_1^* \otimes E_1 \to 0.
\]
Now $\mbox{End}'E \subseteq K$ and the image of $\mbox{End}'E$ in $E_1^* \otimes E_1$ is 
\[
\bigoplus_{i=1}^q \mbox{End}(\CC^{m_{i1}}) \otimes \mbox{End}(D_i),
\]
so since $\mbox{End}''E = \mbox{End}E / \mbox{End}'E$ we get short exact sequences
\begin{equation}
\label{o1}
0 \to \frac{K}{\mbox{End}'E} \to \mbox{End}''E \to E^*_1 \otimes \left(\frac{E}{E_1} \right) \to 0
\end{equation}
and
\begin{equation}
\label{o2}
0 \to \frac{\mbox{End}'E + ((E/E_1)^* \otimes E)}{\mbox{End}'E} \to \frac{K}{\mbox{End}'E} 
\to \frac{E_1^* \otimes E_1}{\bigoplus_{i=1}^q \mbox{End}(\CC^{m_{i1}}) \otimes \mbox{End}(D_i)} \to 0.
\end{equation}
Since
\begin{eqnarray*}
\frac{\mbox{End}'E + ((E/E_1)^* \otimes E)}{\mbox{End}'E} & \cong & \frac{(E/E_1)^* 
\otimes E}{\mbox{End}'E \cap ((E/E_1)^* \otimes E)}\\
& \cong & \frac{((E/E_1)^* \otimes E)/((E/E_1)^* \otimes E_1)}{(\mbox{End}'E \cap ((E/E_1)^* 
\otimes E))/((E/E_1)^* \otimes E_1)}\\
& \cong & \frac{(E/E_1)^* \otimes (E/E_1)}{\mbox{End}'(E/E_1)} = \mbox{End}''(E/E_1),
\end{eqnarray*}
the short exact sequence (\ref{o2}) becomes 
\begin{equation}
\label{o3}
0 \to \mbox{End}''\left( \frac{E}{E_1} \right) \to \frac{K}{\mbox{End}'E} \to \frac{E_1^* 
\otimes E_1}{\bigoplus_{i=1}^q \CC^{m_{i1}^2} \otimes D_i^* \otimes D_i} \to 0.
\end{equation}
From the sequences (\ref{o1}) and (\ref{o3}) it follows that the rank of $\mbox{End}''E$ is
\begin{eqnarray*}
\mbox{rank}(\mbox{End}''E) & = & \mbox{rank}(K/\mbox{End}'E) + \mbox{rank}(E_1^* \otimes (E/E_1))\\
& = & \mbox{rank}(\mbox{End}''(E/E_1)) + \mbox{rank}(E_1^* \otimes E_1) - \sum_{i=1}^q (m_{i1})^2 
\mbox{rank}(D_i^* \otimes D_i)\\
&  &  + \mbox{ rank}(E_1^* \otimes (E/E_1))\\
& = & \mbox{rank}(\mbox{End}''(E/E_1)) + \mbox{rank}(E_1^* \otimes E)  - \sum_{i=1}^q (m_{i1})^2 (n_i)^2.
\end{eqnarray*}
Thus by induction on $r$ we have 
\[
\mbox{rank}(\mbox{End}''E) = \sum_{i_1,i_2=1}^{q} \sum_{1 \leq j_1\leq j_2\leq r} m_{i_1 j_1} m_{i_2 j_2} 
n_{i_1} n_{i_2} - \sum_{i=1}^{q} \sum_{j=1}^{r} (n_i)^2.
\]
Since $D_1,\ldots,D_q$ all have the same slope $\D/\n$ as $E$ we have 
$\mbox{deg}(\mbox{End}''E) =0.$ Therefore by Riemann-Roch
$$ \dim H^1 (\Sigma, \mbox{End}''E) = \dim H^0 (\Sigma,\mbox{End}''E) $$
$$ + (g-1)\left(  \sum_{i_1,i_2=1}^{q} \sum_{1 \leq j_1\leq j_2\leq r} m_{i_1 j_1} m_{i_2 j_2} 
n_{i_1} n_{i_2} - \sum_{i=1}^{q} \sum_{j=1}^{r} (n_i)^2
\right) .   $$
Moreover the short exact sequences (\ref{o1}) and (\ref{o3}) give us long exact sequences of cohomology
\[
0 \to H^0(\Sigma, K/\mbox{End}'E)  \to H^0(\Sigma,\mbox{End}''E) \to H^0(\Sigma, E_1^* \otimes (E/E_1)) \to \cdots
\]
and
$$0 \to H^0(\Sigma,\mbox{End}''(E/E_1)) \to H^0(\Sigma, K/\mbox{End}'E) $$
$$\to H^0(\Sigma, 
E_1^* \otimes E_1/\bigoplus_{i=1}^q \mbox{End}(\CC^{m_{ij}}) \otimes \mbox{End}(D_i)) \to \cdots
$$
Now $E_1 \cong \bigoplus_{i=1}^q \CC^{m_{i1}} \otimes D_i$ where $D_1,...,D_q$ are
nonisomorphic stable bundles all of the same slope as $E_1$, and
$$E_2/E_1 \cong \bigoplus_{i=1}^q \CC^{m_{i2}} \otimes D_i$$
is the maximal subbundle of $E/E_1$ which is a direct sum of stable bundles all of the same
slope as $E/E_1$. Since $E_1$ and $E/E_1$ have the same slope, it follows from Corollary \ref{lem34}
that
$$H^0(\Sigma, E_1^* \otimes (E/E_1)) = H^0(\Sigma, E_1^* \otimes (E_2/E_1)) \cong \bigoplus_{i=1}^q
(\CC^{m_{i1}})^* \otimes \CC^{m_{i2}}.$$
By choosing an open cover $\mathcal{U}$ of $\Sigma$ such that the filtration 
$0 = E_0 \subset E_1 \subset E_2 \subset \cdots \subset E_r = E$
is trivial over each $U \in \mathcal{U}$, and describing $E$ in terms of upper triangular
transition functions on $E_1 \oplus (E_2/E_1) \oplus ... \oplus (E/E_{r-1})|_{U\cap V}$ for
$U,V \in \mathcal{U}$ which induce the identity on $E_j/E_{j-1}$ for $1 \leq j \leq r$, we see
that the natural map
$$H^0(\Sigma, \mbox{End}''E) \to H^0(\Sigma, E_1^* \otimes (E/E_1)) = H^0(\Sigma, E_1^* \otimes (E_2/E_1))$$
is surjective. Also
$$H^0(\Sigma, E_1^* \otimes E_1 / \bigoplus_{i=1}^q \mbox{End}(\CC^{m_{i1}} \otimes \mbox{End}(D_i))
= \bigoplus_{i \neq j} (\CC^{m_{i1}})^* \otimes \CC^{m_{i1}} \otimes H^0(\Sigma, D_i^* \otimes D_j) = 0,$$
so
$$\dim H^0(\Sigma, \mbox{End}''E) = \sum_{i=1}^q m_{i1}m_{i2} + \dim H^0(\Sigma, \mbox{End}''(E/E_1))$$
and thus by induction on $r$ we have
$$\dim H^0(\Sigma, \mbox{End}''E) = \sum_{i=1}^q \sum_{j=1}^{r-1} m_{ij} m_{i j+1}.$$
Therefore $
\dim H^1(\Sigma, \mbox{End}''E) $ is equal to
\begin{equation} \label{dimeq} \sum_{i=1}^q \sum_{j=1}^{r-1} m_{ij} m_{i j+1} + (g-1) \left( 
\sum_{i_1,i_2 = 1}^q \sum_{1\leq j_1 \leq j_2 \leq r} m_{i_1 j_1 } m_{i_2 j_2} n_{i_1} n_{i_2}
- \sum_{i=1}^q \sum_{j=1}^r (m_{ij} n_i)^2
\right). \end{equation}
In particular this tells us that the image in $\mathcal{C}$ of the derivative of the
injection (\ref{inj}) has constant codimension, and it follows as in \cite[$\S$7]{AB} 
(see also \cite[$\S \S$14 and 15]{AB}) that the subset $S_{[{\bf n,m}]}$ is locally a complex submanifold
of $\mathcal{C}$ of finite codimension given by (\ref{dimeq}).

\begin{cor} \label{4.6} If $[{\bf d,n,m}] \in \mathcal{I}$ then $S_{[{\bf d,n,m}]}$ is a locally closed complex
submanifold of $\mathcal{C}$ of codimension
$$ \sum_{1 \leq k_2 < k_1 \leq s} (n_{k_1} d_{k_2} - n_{k_2} d_{k_1} + n_{k_1} n_{k_2}(g-1)) +
\sum_{k=1}^s \sum_{i=1}^{q}\sum_{j=1}^{r-1} m_{ijk} m_{i j+1 k}$$
$$ + (g-1) \sum_{k=1}^s \left( \sum_{i_1,i_2=1}^{q} \sum_{1 \leq j_1 \leq j_2\leq r} 
m_{i_1 j_1 k} m_{i_2 j_2 k} n_{i_1 k} n_{i_2 k} - \sum_{i=1}^{q} \sum_{j=1}^{r} (m_{ijk}n_{ik})^2 \right) .
$$
\end{cor}

\noindent{\bf Proof}: This follows immediately from (\ref{12}), Lemma \ref{4.4} and the 
definition of $S_{[{\bf d,n,m}]}$ (Definition \ref{def81}).

\begin{rem} \label{4.7} Let
$0=E_0 \subset E_1 \subset ... \subset E_t = E$ be the maximal Jordan--H\"{o}lder
filtration of a bundle $E$. Then the kernels of the duals of the inclusions
$E_j \to E$ give us a filtration
$$0=F_0 \subset F_1 \subset ... \subset F_t = E'$$
of the dual $E'$ of $E$, such that if $1 \leq j \leq t$ then
$$F_j/F_{j-1} \cong (E_{t-j+1}/E_{t-j})'.$$
Thus $F_j/F_{j-1}$ is a direct sum of stable bundles all of the same slope, and moreover $F_{j-1}$ is
the minimal subbundle of $F_j$ such that $F_j/F_{j-1}$ is a direct sum of stable bundles all
of which have minimal slope among quotients of $F_j$. 

Applying this construction with $E$ replaced by $E'$, we find that every holomorphic
bundle $E$ over $\Sigma$ has a canonical filtration
$$0=F_0 \subset F_1 \subset ... \subset F_t = E,$$
which we will call the {\em minimal Jordan--H\"{o}lder filtration} of $E$, with
the property that if $1 \leq j \leq t$ then $F_{j-1}$ is the minimal
subbundle of $F_j$ such that $F_j /F_{j-1}$ is a direct sum of stable
bundles all of which have minimal slope among quotients of $F_j$.

The minimal and maximal Jordan--H\"{o}lder filtrations of a bundle
do not necessarily coincide. For example, consider the direct sum $E \oplus F$
of two semistable bundles with maximal Jordan--H\"{o}lder filtrations
$0=E_0 \subset E_1 \subset ... \subset E_t = E$
and
$0=F_0 \subset F_1 \subset ... \subset F_s = F$
where without loss of generality we may assume that $s \leq t$. 
If $E$ and $F$ have the same slope, then it is easy to check that
the maximal Jordan--H\"{o}lder filtration of $E \oplus F$ is
$$0=E_0 \oplus F_0 \subset E_1 \oplus F_1 \subset ... \subset E_s \oplus F_s \subset
E_{s+1} \oplus F_s \subset ... \subset E_t \oplus F_s;$$
that is, it is the direct sum of the maximal Jordan--H\"{o}lder filtrations
of $E$ and $F$ with the shorter one extended trivially at the top. Similarly the
minimal Jordan--H\"{o}lder filtration of $E \oplus F$ is
the direct sum of the minimal Jordan--H\"{o}lder filtrations
of $E$ and $F$ with the shorter one extended trivially at the bottom. Thus if
the minimal and maximal Jordan--H\"{o}lder filtrations of $E$ and $F$ coincide
(which will be the case if, for example, each of the subquotients $E_j/E_{j-1}$ and
$F_j/F_{j-1}$ are stable) but these filtrations are not of the same length, then the
minimal Jordan--H\"{o}lder filtration
$$0=E_0\oplus F_0 \subset E_1 \oplus F_0 \subset ... \subset E_{t-s} \oplus F_0
\subset E_{t-s+1} \oplus F_1 \subset ... \subset E_t \oplus F_s$$
of $E \oplus  F$ will be different from its maximal Jordan--H\"{o}lder filtration.
\end{rem}

\begin{df} \label{4.8} Given $[{\bf d,n,m}] \in \mathcal{I}$, let $S^{minJH}_{[{\bf d,n,m}]}$ denote
the subset of $\mathcal{C}$ consisting of those holomorphic 
structures on our fixed smooth bundle of rank $\n$ and degree $\D$ whose minimal Jordan--H\"{o}lder
filtration is of the form
\[
0 = E_{0,1} \subset E_{1,1} \subset \cdots \subset E_{r_1,1}=E_{0,2}\subset E_{1,2} \subset \cdots \]
\[ \cdots
\subset E_{r_{s-1}, s-1} = E_{r_s, 0} \subset \cdots \subset E_{r_s,s} = E
\]
with
\[
E_{j,k}/E_{j-1,k} \cong (\CC^{m_{1jk}} \otimes D_{1k}) \oplus \cdots \oplus (\CC^{m_{q_k jk}} \otimes D_{q_k k})
\] 
for $1 \leq k \leq s$ and $1 \leq j \leq q_k$, where $D_{1 k},\ldots,D_{q_k k}$ are nonisomorphic stable bundles with 
\[
\mbox{rank}(D_{ik}) = n_{ik} \mbox{ and } \mbox{deg}(D_{ik}) = d_{ik}
\]
and $E_{j-1,k}$ is the minimal subbundle of $E_{j,k}$ such that $E_{j,k}/E_{j-1,k}$ is a direct sum of stable 
bundles of  slope $d_k/n_k$. 
\end{df}

\section{More indexing sets}
\renorm

 In this section we will  consider the
indexing set $\Gamma$ for the stratification $\{\Sigma_\gamma: \gamma \in \Gamma\}$
of $\mathcal{C}^{ss}$ defined as in $\S$2. 

If $\gamma \in \Gamma$ then by (\ref{new*}) either $\gamma = 0$ or $\gamma = R_l$ 
or $\gamma \in \mathcal{B}_l \backslash \{0\}\times \{ l \}$ for some $1 \leq l \leq \tau$. 
If $\gamma = 0$ then $\Sigma_\gamma = \mathcal{C}^s$, while by 
\cite[pp.248-9]{K5} if $\gamma=R_l$ then there exists
$ [{\bf n},{\bf m}] = [(n_i)_{i=1}^{q},(m_{ij})_{i=1,j=1}^{q,r}] \in \mathcal{I}^{ss} $
with $r=1$ and $q=q_1$, such that 
\begin{equation} \label{rl}
R_l = \prod_{i=1}^{q} GL(m_{i};\CC)
\end{equation}
where $m_i = m_{i1}$, and $\Sigma_{R_l}$ consists of all those holomorphic structures $E$ with
\begin{equation}
\label{newo**}
E \cong (\CC^{m_1} \otimes D_1) \oplus \cdots \oplus (\CC^{m_q} \otimes D_q)
\end{equation}
where $D_1, \ldots ,D_q$ are all stable of slope $\D/\n$ and not isomorphic to one another. 

In order to describe the strata $\{\Sigma_{\beta,l}: \beta \in \mathcal{B}_l \backslash \{0\}\}$ more 
explicitly, we need to look at the action of $R_l$ on the normal $\mathcal{N}_{R_l}$ to $GZ^{ss}_{R_l}$ at a point 
represented by a holomorphic structure $E$ of the form (\ref{newo**}), 
and to understand the stratification on $\PP(\mathcal{N}_{R_l})$
induced by this action of $R_l$. If we choose a $C^{\infty}$ isomorphism of
our fixed $C^{\infty}$ bundle $\mathcal{E}$ with $ (\CC^{m_1} \otimes D_1) \oplus \cdots \oplus (\CC^{m_q} \otimes D_q)$,
then we can identify $\mathcal{C}$ with the infinite-dimensional vector space
$$\Omega^{0,1}(\mbox{End}( (\CC^{m_1} \otimes D_1) \oplus \cdots \oplus (\CC^{m_q} \otimes D_q)))$$
and the normal to the $\mathcal{G}_c$-orbit at $E$ can be identified with $H^1(\Sigma, \mbox{End} E)$,
where $\mbox{End} E$ is the bundle of holomorphic endomorphisms of $E$ \cite[$\S$7]{AB}.
If $\delta_i^j$ denotes the Kronecker delta then 
the normal to $G Z_R^{ss}$ can be identified with
\[
H^1(\Sigma,\mbox{End}'_{\oplus} E) \cong \bigoplus_{i_1,i_2=1}^{q} 
\CC^{m_{i_1}m_{i_2}-\delta_{i_1}^{ i_2}} \otimes H^1(\Sigma,D^*_{i_1} \otimes D_{i_2})
\]
where $\mbox{End}'_{\oplus} E$ is the quotient of the bundle $\mbox{End}E$ 
of holomorphic endomorphisms of $E$ by the subbundle $\mbox{End}_{\oplus} E$ 
consisting of those endomorphisms which preserve the decomposition (\ref{newo**}). The action of 
$ R_l = \prod_{i=1}^{q} GL(m_i;\CC) $
on this is given by the natural action on $\CC^{m_{i_1} m_{i_2} - \delta_{i_1}^{ i_2}}$ identified 
with the set of $m_{i_1} \times m_{i_2}$ matrices if $i_1 \neq i_2$ and the set of trace-free $m_{i_1} \times m_{i_1}$
matrices if $i_1 = i_2$; its weights $\alpha$ are therefore of the form $\alpha = \xi - \xi'$ where $\xi$ 
and  $\xi'$ are weights of the standard representation of $R_l$ on $\oplus_{i=1}^{q} \CC^{m_i}$ (see 
\cite[pp.251-2]{K5} noting the error immediately before (3.18)).

Any element $\beta$ of the indexing set $\mathcal{B}_l$ 
 is represented by the closest point to 0 of the convex hull of some nonempty set of these weights,
and two such closest points can be taken to represent the same element of $\mathcal{B}_l$ if and only if they lie in the
same $Ad(N_l)$-orbit, where $N_l$ is the normaliser of $R_l$ (see \cite{K2} or \cite{K7}). 
By (\ref{N_0}) the orbit of $\beta$ under
the adjoint action of the connected component of $N_l$ is just its $Ad(R_l)$-orbit, and so by
(\ref{N/N_0}) the $Ad(N_l)$-orbit of $\beta$ is the union
$$\bigcup_{w \in \pi_0(N_l)} w. Ad(N_l) (\beta)$$
where $\pi_0(N_l)$ is the product of permutation groups
$$\pi_0(N_l) = \prod_{j \geq 0, k \geq 0} Sym(\#\{i:m_i = j \mbox{ and } n_i=k\}).$$

We can describe this indexing set $\mathcal{B}_l$ more explicitly as follows. Let us take our maximal compact torus
$T_l$ in $R_l$ to be the product of the standard maximal tori of the unitary groups $U(m_1)$,..., $U(m_q)$
consisting of the diagonal matrices, and let $\liet_l$ be its Lie algebra. Let
$$M = m_1 + ... + m_q$$
and let $e_1,...,e_M$ be the weights of the standard representation of $T_l$ on
$\CC^{m_1}\oplus ... \oplus \CC^{m_q}$. We take the usual invariant inner product on
the Lie algebra ${\bf u}(p)$ of $U(p)$ given by $\langle A, B\rangle = -{\rm tr} A\bar{B}^t$
and restrict it to $T_l$. Since $R_l$ is embedded in $GL(p;\CC)$ by identifying
$\oplus_{i=1}^q (\CC^{m_i} \otimes \CC^{p_i})$ with $\CC^p$, it follows that $e_1,...,e_M$
are mutually orthogonal and $|\!|e_j|\!|^2 = 1/p_i$ if $m_1 +...+ m_{i-1} < j \leq m_1 +...+m_i$.

\begin{prop} \label{5.1} Let $\beta$ be any nonzero element of the Lie algebra
$\liet_l$ of the maximal compact torus $T_l$ of $R_l$. Then $\beta$ represents
an element of $\mathcal{B}_l \backslash \{ 0 \}$ if and only if there is a partition
$$\{ \Delta_{h,m}:(h,m) \in J \}$$
of $\{1,...,M \}$, indexed by a subset $J$ of $\ZZ \times \ZZ$
of the form
$$J = \{(h,m) \in \ZZ \times \ZZ: 1 \leq h \leq L, l_1(h) \leq m \leq l_2(h) \}$$
for some positive integer $L$ and functions $l_1$ and $l_2:\{1,...,L\} \to \ZZ$
such that $l_1(h) \leq l_2(h)$ for all $h \in \{1,...,L\}$, with the following properties.
If
$$r_{h,m} = \sum_{j \in \Delta_{h,m}} |\!|e_j|\!|^{-2}$$
then the function $\epsilon: \{ 1,...,L\} \to \QQ$ defined by
$$\epsilon(h) = \left( \sum_{m=l_1(h)}^{l_2(h)} m r_{h,m}\right) 
\left(\sum_{m=l_1(h)}^{l_2(h)} r_{h,m}\right)^{-1}$$
satisfies $-1/2 \leq \epsilon(h) < 1/2$ and $\epsilon(1) > \epsilon(2) >...> \epsilon(L)$, and
$$\frac{\beta}{|\!| \beta |\!|^2} = \sum_{h=1}^L \sum_{m = l_1(h)}^{l_2(h)}
\sum_{j \in \Delta_{h,m}} (\epsilon(h) - m) \frac{e_j}{|\!| e_j |\!|^2}.$$
\end{prop}

\begin{rem} \label{9.2} Note that because of the conditions on the function $\epsilon$,
the partition $\{\Delta_{h,m}:(h,m) \in J \}$ and its indexing can be recovered from
the coefficients of $\beta$ with respect to the basis $e_1/|\!|e_1 |\!|^2, ..., e_M/|\!|e_M|\!|^2$
of $\liet_l$.
\end{rem}

\newcommand{\lij}{\lambda_{ij}}
\newcommand{\lji}{\lambda_{ji}}
\newcommand{\lik}{\lambda_{ik}}
\newcommand{\ljk}{\lambda_{jk}}
\newcommand{\lki}{\lambda_{ki}}
\newcommand{\lkj}{\lambda_{kj}}
\newcommand{\lijbh}{\lambda_{ij}^{\beta h}}
\newcommand{\lijb}{\lambda_{ij}^{\beta}}

\noindent {\bf Proof of Proposition \ref{5.1}}:  $\beta \in \liet_l  \backslash \{ 0 \}$ represents
an element of $\mathcal{B}_l \backslash \{ 0\}$ if and only if it is the closest point to 0
of the convex hull of
$$\{ e_i - e_j: (i,j) \in S \}$$
for some nonempty subset $S$ of $\{ (i,j) \in \ZZ \times \ZZ: 1 \leq i,j \leq M \}$. Then 
$\beta$ can be expressed in the form
$$ \beta = \sum_{(i,j) \in S} \lambda_{ij}^{\beta} (e_i - e_j)$$
for some $\lambda_{ij}^{\beta} \in \RR$ for $(i,j) \in S$ such that $\lambda_{ij}^{\beta} \geq 0$ and
$\sum_{(i,j) \in S} \lambda_{ij}^{\beta} =1$. Replacing $S$ with its subset
$\{ (i,j) \in S: \lijb > 0 \}$ we may assume that $\lijb > 0$ for all $(i,j) \in S$.
Moreover clearly if $S = \{ e_i - e_j \}$ has just one element then $\beta = e_i - e_j$,
and we then take $J = \{ (1,0), (1,1) \}$ with $\Delta_{(1,0)} = \{ j\}$ and
$\Delta_{(1,1)} = \{ i \}$, so we can assume without loss of generality that
$\lijb < 1$ for all $(i,j) \in S$.  Since $\beta \neq 0$ we can also assume that
the convex hull of $\{ e_i - e_j: (i,j) \in S \}$ does not contain 0.

In order to find the closest point to 0 of the convex hull of $\{ e_i - e_j: (i,j) \in S \}$
we minimise
$$|\!| \sum_{(i,j) \in S} \lij (e_i - e_j) |\!|^2$$
subject to the constraints that $\lij \geq 0$ for all $(i,j) \in S$ and $\sum_{(i,j)\in S}
\lij = 1$. Since the weights $e_1,...,e_M$ are mutually orthogonal, we have
$$|\!| \sum_{(i,j) \in S} \lij (e_i - e_j) |\!|^2 = |\!| \sum_{i=1}^M (\sum_{j:(i,j) \in S} \lij 
- \sum_{j:(j,i) \in S} \lji) e_i |\!|^2 $$
$$= \sum_{i=1}^M (\sum_{j:(i,j) \in S} \lij 
- \sum_{j:(j,i) \in S} \lji)^2 |\!| e_i |\!|^2.$$
Using the method of Lagrange multipliers, we consider
$$ \sum_{i=1}^M (\sum_{j:(i,j) \in S} \lij 
- \sum_{j:(j,i) \in S} \lji)^2 |\!| e_i |\!|^2 - \lambda(\sum_{(i,j) \in S} \lij - 1).$$
If $(i, j) \in S$ then $i \neq j$ and $(j, i) \not \in S$ since the
convex hull of $\{ e_i - e_j: (i,j) \in S \}$ does not contain 0, so
$$\frac{\partial}{\partial \lambda_{i j}} (\sum_{i=1}^M (\sum_{j:(i,j) \in S} \lij 
- \sum_{j:(j,i) \in S} \lij)^2 |\!| e_i |\!|^2 - \lambda(\sum_{(i,j) \in S} \lij - 1) )$$
is equal to 
$$2(\sum_{k:(i,k) \in S} \lik
- \sum_{k:(k,i) \in S} \lki) |\!| e_i |\!|^2 - 2(\sum_{k:(j,k) \in S} \ljk - \sum_{k:(k,j) \in S}
\lkj )|\!| e_j|\!|^2 - \lambda.$$
Thus $\beta = \sum_{(i,j) \in S} \lijb (e_i - e_j)$ where for each $(i,j) \in S$ we have
either $\lijb = 0$ or $\lijb=1$ (both of which are ruled out by the assumptions on $S$)
or
\begin{equation} \label{star} (\sum_{k:(i,k) \in S} \lik
- \sum_{k:(k,i) \in S} \lki) |\!| e_i |\!|^2 - (\sum_{k:(j,k) \in S} \ljk - \sum_{k:(k,j) \in S}
\lkj )|\!| e_j|\!|^2 = \lambda/2 \end{equation}
where $\lambda$ is independent of $(i,j) \in S$.

From $S$ we can construct a directed graph $G(S)$ with vertices $1,...,M$ and directed
edges from $i$ to $j$ whenever $(i,j)\in S$. 
Let $\Delta_1$,..., $\Delta_N$ be the connected components of this graph.
Then $\{ e_i - e_j: (i,j) \in S \}$ is the disjoint union of its subsets $\{ e_i - e_j: (i,j) \in S, i,j \in \Delta_h
 \}$ for $1 \leq h \leq N$, and $\{ e_i - e_j: (i,j) \in S, i,j \in \Delta_h
 \}$ is contained in the vector subspace of $\liet_l$ spanned by the basis
vectors $\{e_k: k \in \Delta_h\}$. Since these subspaces are mutually
orthogonal for $1 \leq h \leq N$, it follows that
\begin{equation} \label{beta} \beta =\left(   \sum_{h=1}^L \frac{1}{|\!| \beta_h |\!|^2}
\right )^{-1} \sum_{h=1}^L \frac{\beta_h}{|\!| \beta_h |\!|^2} 
\end{equation}
where 
$$\beta_h = \sum_{(i,j) \in S, i,j \in\Delta_h} \lijbh (e_i - e_j)$$
is the closest point to 0 of the convex hull of
$\{ e_i - e_j: (i,j) \in S, i,j \in \Delta_h  \}$ for $1 \leq h \leq N$, and
without loss of generality we assume that $\beta_h$ is nonzero
when $1 \leq h \leq L$ and zero when $L < h \leq N$. Note that then
\begin{equation} \label{betasq} |\!|\beta|\!|^2  =\left(   \sum_{h=1}^L \frac{1}{|\!| \beta_h |\!|^2}
\right )^{-2} \sum_{h=1}^L \frac{|\!|\beta_h|\!|^2}{|\!| \beta_h |\!|^4} 
 =\left(   \sum_{h=1}^L \frac{1}{|\!| \beta_h |\!|^2}
\right )^{-1}
\end{equation}
so that
\begin{equation} \label{beta2} \frac{\beta}{|\!|\beta|\!|^2} = \sum_{h=1}^L \frac{\beta_h}{|\!| \beta_h |\!|^2} 
\end{equation}
and
$$\lijbh = (|\!| \beta_h|\!|^2 / |\!| \beta |\!|^2) \lijb$$
if $i, j \in \Delta_h$.
 For
$1 \leq h \leq L$ let $\kappa_h$ be defined by
$$\kappa_h = {\rm max} \{ (\sum_{k:(i,k) \in S} \lik^{\beta}
- \sum_{k:(k,i) \in S} \lki^{\beta}) |\!| e_i |\!|^2: i \in \Delta_h \}$$
where a sum over the empty set is interpreted as 0. Then by (\ref{star}) for
$1 \leq h \leq L$ we can express $\Delta_h$ as a disjoint union
$$\Delta_h = \hat{\Delta}_{h,0} \sqcup  \hat{\Delta}_{h,1} \sqcup \ldots \sqcup \hat{\Delta}_{h,l_h} $$
where
\begin{equation} \label{dagger} \hat{\Delta}_{h,m} = \{ i \in \Delta_h : (\sum_{j:(i,j) \in S} \lij^{\beta}
- \sum_{j:(j,i) \in S} \lji^{\beta}) |\!| e_i |\!|^2 = \kappa_h - m|\lambda|/2\}. \end{equation}
Let us assume that $\lambda \leq 0$; the argument is similar if $\lambda \geq 0$. Then (\ref{star})
tells us that
if $(i,j)\in S$  then there exist $h$ and $m$ such that $i \in \hat{\Delta}_{h,m}$
and $j \in \hat{\Delta}_{h,m+1}$.
Note that if $\hat{\Delta}_{h,m_1}$ and $\hat{\Delta}_{h,m_2}$ are
nonempty then so is $\hat{\Delta}_{h,m}$ whenever $m_1 < m < m_2$, so without
loss of generality we may assume that $\hat{\Delta}_{h,m}$ is nonempty when
$1 \leq h \leq L$ and $0 \leq m \leq l_h$.
Thus if $1 \leq h \leq L$ we have
$$\beta_h = (|\!| \beta_h|\!|^2 / |\!| \beta |\!|^2)
\sum_{i \in \Delta_h} ( \sum_{j \in \Delta_h, (i,j) \in S} \lijb
- \sum_{j \in \Delta_h, (j,i) \in S} \lji^{\beta}) e_i$$
$$= (|\!| \beta_h|\!|^2 / |\!| \beta |\!|^2)
\sum_{m=1}^{l_h} \sum_{i \in \hat{\Delta}_{h,m}} (\kappa_h - m|\lambda|/2) \frac{e_i}{|\!|e_i |\!|^2}.$$
For $1 \leq h \leq L$ and $0 \leq m \leq l_h$ let
$r_{h,m} = \sum_{j \in \hat{\Delta}_{h,m}} |\!| e_j |\!|^{-2};$
then by (\ref{dagger}) we have
$$ \sum_{i \in \hat{\Delta}_{h,m}} ( \sum_{j:(i,j) \in S} \lij^{\beta}
- \sum_{j:(j,i) \in S} \lji^{\beta})  = r_{h,m} (\kappa_h - m|\lambda|/2\}).$$
Recall that if $(i,j) \in S$ then $i \in \hat{\Delta}_{h,m}$ if and only
if $j  \in \hat{\Delta}_{h,m+1}$, so we get
$$ \sum_{i \in \hat{\Delta}_{h,0}}  \sum_{j:(i,j) \in S} \lij^{\beta}
 = r_{h,0} \kappa_h,$$
and hence
$$ \sum_{i \in \hat{\Delta}_{h,1}}  \sum_{j:(i,j) \in S} \lij^{\beta} = 
\sum_{i \in \hat{\Delta}_{h,1}} ( \sum_{j:(i,j) \in S} \lij^{\beta}
- \sum_{j:(j,i) \in S} \lji^{\beta}) + \sum_{j \in \hat{\Delta}_{h,0}}\sum_{j:(j,i) \in S} \lji^{\beta}$$
$$ = r_{h,1} (\kappa_h - |\lambda|/2\}) + r_{h,0} \kappa_h,$$
and similarly
$$\sum_{i \in \hat{\Delta}_{h,m}}  \sum_{j:(i,j) \in S} \lij^{\beta} = 
(r_{h,0} + r_{h,1} + ... + r_{h,m}) \kappa_h - (r_{h,1} + 2r_{h,2} + ...
+ mr_{h,m})\frac{|\lambda|}{2}$$
if $1 \leq m \leq l_h$. Since $ \sum_{(i,j) \in S, i,j \in\Delta_h}  \lijbh = 1$, 
it follows that
$$|\!| \beta|\!|^2 / |\!| \beta_h |\!|^2 = ((l_h+1)r_{h,0} + l_h r_{h,1} +...+ r_{h,l_h})
\kappa_h $$
$$- (l_h r_{h,1} + 2(l_h-1)r_{h,2} + 3(l_h - 2) r_{h,3} + ... + l_h r_{h,l_h})|\lambda|/2,$$
and since  $ \sum_{(i,j) \in S, i,j \in\Delta_h}  \lijb -   \sum_{(j,i) \in S, i,j \in\Delta_h}  \lji^{\beta} = 0$ we have
$$0 = (r_{h,0}+r_{h,1} + ... + r_{h,l_h}) \kappa_h - (r_{h,1} + 2r_{h,2} + ... + l_h r_{h,l_h})|\lambda|/2.$$
Thus
$$|\lambda|/2 = (|\!| \beta|\!|^2 / |\!| \beta_h |\!|^2) (r_{h,0} + r_{h,1} + ... + r_{h,l_h})/\mu_h$$
and
$$\kappa_h = (|\!| \beta|\!|^2 / |\!| \beta_h |\!|^2) (r_{h,1} + 2 r_{h,2} + ... + l_h r_{h,l_h})/\mu_h$$
where
$$\mu_h = \sum_{i,j=0}^{l_h} ((l_h - i +1)j - j(l_h - j + 1)) r_{h,i} r_{h,j}$$
$$= \sum_{0\leq i<j\leq l_h} ((j-i)j + (i-j)i) r_{h,i} r_{h,j}
= \sum_{0\leq i<j\leq l_h} (j-i)^2 r_{h,i} r_{h,j}.$$
Therefore 
$$\beta_h = \sum_{m=0}^{l_h} \sum_{j=0}^{l_h} \frac{(j-m)r_{h,j}}{\mu_h} 
\sum_{i \in \hat{\Delta}_{h,m}} \frac{e_i}{|\!| e_i |\!|^2}$$
and 
$$|\!| \beta_h |\!|^2 = (r_{h,0} + r_{h,1} +...
+ r_{h,l_h}))/\mu_h .$$
By defining 
$\Delta_{h,m} = \hat{\Delta}_{h,m - l_1(h)}$ for an
appropriate integer $l_1(h)$, we can arrange that the function $\epsilon$
defined in the statement of the proposition takes values in the interval $[-1/2,1/2)$,
and then by amalgamating those ${\Delta}_{h}$ for which $\epsilon(h)$ takes the same value
and rearranging them so that $\epsilon$ is a strictly decreasing function, we
can assume that the required conditions on $\epsilon$ are satisfied, and we have
$$\frac{\beta}{|\!| \beta |\!|^2} = \sum_{h=1}^L \sum_{m = l_1(h)}^{l_2(h)}
\sum_{i \in \Delta_{h,m}} (\epsilon(h) - m) \frac{e_i}{|\!| e_i |\!|^2}.$$
This gives us all the required properties if $L=N$; that is, if $\bigcup_{(h,m) \in J} \Delta_{h,m}$ is equal to 
$\{1,...,M\}$. Otherwise we amalgamate $\{1,...,M\} \backslash \bigcup_{(h,m) \in J} \Delta_{h,m}$
with $\Delta_{h_0,m_0}$ where $(h_0,m_0)$ is the unique element of $J$ such that
$\epsilon(h_0) = 0 = m_0$ if such an element exists, and if there is no such element of
$J$ then we adjoin $(L+1,0)$ to $J$ and define
$$\Delta_{L+1,0} = \{1,...,M\} \backslash \bigcup_{h=1}^L \bigcup_{m=l_1(h)}^{l_2(h)} \Delta_{h,m}.$$

Conversely, suppose that we are given any partition $\{ \Delta_{h,m}: (h,m) \in J\}$ of $\{1,...,M\}$
indexed by
$$J = \{ (h,m) \in \ZZ \times \ZZ: 1 \leq h \leq L, l_1(h) \leq m \leq l_2(h) \}$$
for some positive integer $L$ and functions $l_1$ and $l_2: \{1,...,L\} \to \ZZ$ with
$l_1 \leq l_2$, satisfying $\epsilon(h) \in [-1/2,1/2)$ and
$\epsilon(1) > \epsilon(2) > ... > \epsilon(L)$ where
$$\epsilon(h) =\left( \sum_{m=l_1(h)}^{l_2(h)} m r_{h,m} \right) 
\left(\sum_{m=l_1(h)}^{l_2(h)} r_{h,m}\right)^{-1}$$
for $r_{h,m} = \sum_{i \in \Delta_{h,m}} |\!|e_i|\!|^{-2}$. Suppose also
that 
$$\beta = \hat{\beta}/ |\!|\hat{\beta}|\!|^2$$
(or equivalently $\hat{\beta} = \beta / |\| \beta |\|^2$) where
$$\hat{\beta}  = \sum_{h=1}^L \sum_{m = l_1(h)}^{l_2(h)}
\sum_{j \in \Delta_{h,m}} (\epsilon(h) - m) \frac{e_j}{|\!| e_j |\!|^2}.$$
It suffices to show that $\beta $ is the closest point to 0 of the convex hull of
$\{e_i - e_j: (i,j) \in S \}$, where $S$ is the set of ordered pairs $(i,j)$ with
$i,j \in \{1,...,M\}$ such that $e_i \in \Delta_{h,m}$ and $e_j \in \Delta_{h,m+1}$
for some $(h,m) \in J$ such that $(h,m+1) \in J$. For this, it is enough to prove
firstly that $\beta$ lies in the convex hull of $\{e_i - e_j: (i,j) \in S \}$ and secondly
that $(e_i - e_j).\beta = |\!|\beta |\!|^2$ (or equivalently that $(e_i - e_j).\hat{\beta} = 1$)
for all $(i,j) \in S$. The latter follows easily from the choice of $S$: if $(i,j) \in S$ then
there exists $(h,m) \in J$ such that $(h,m+1) \in J$ and
$e_i \in \Delta_{h,m}$ and $e_j \in \Delta_{h,m+1}$, so
$$(e_i - e_j ). \hat{\beta} = \epsilon(h) - m - \epsilon(h) +m + 1 = 1.$$
To show  that $\beta$ lies in the convex hull of $\{e_i - e_j: (i,j) \in S \}$, we
note that
$$\sum_{m=l_1(h)}^{l_2(h) -1} \sum_{i \in \Delta_{h,m}} \sum_{j \in \Delta_{h,m+1}}
\sum_{k=l_1(h)}^m \frac{\epsilon(h) r_{h,k} - k r_{h,k}}{r_{h,m} r_{h,m+1} |\!|e_i|\!|^2 |\!|e_j|\!|^2}
(e_i - e_j) $$
$$ = \sum_{m=l_1(h)}^{l_2(h) } \sum_{j \in \Delta_{h,m}}(\epsilon(h) - m) \frac{e_j}{|\!|e_j|\!|^2}.$$
This means that
$$\beta = \sum_{(i,j) \in S} \lijb (e_i - e_j)$$
where
\begin{equation} \label{lam}
\frac{\lijb}{|\!| \beta |\!|^2} = 
\sum_{k=l_1(h)}^m \frac{\epsilon(h) r_{h,k} - k r_{h,k}}{r_{h,m} r_{h,m+1} |\!|e_i|\!|^2 |\!|e_j|\!|^2}
\end{equation}
if $i \in \Delta_{h,m}$ and $j \in \Delta_{h,m+1}$ for some $(h,m) \in J$ such that
$(h,m+1) \in J$. Then
$$\sum_{(i,j) \in S} 
\frac{\lijb}{|\!| \beta |\!|^2} = \sum_{h=1}^L 
\sum_{m=l_1(h)}^{l_2(h) -1} \sum_{i \in \Delta_{h,m}} \sum_{j \in \Delta_{h,m+1}}
\sum_{k=l_1(h)}^m \frac{\epsilon(h) r_{h,k} - k r_{h,k}}{r_{h,m} r_{h,m+1} |\!|e_i|\!|^2 |\!|e_j|\!|^2}$$
$$= \sum_{h=1}^L \sum_{k=l_1(h)}^{l_2(h) -1} \sum_{m=k}^{l_2(h)-1} (\epsilon(h) - k)r_{h,k}
 = \sum_{h=1}^L \sum_{k=l_1(h)}^{l_2(h) -1} (l_2(h)-k) (\epsilon(h) - k)r_{h,k}.$$
Expanding the brackets, using the definition of $\epsilon$ and replacing the index $k$ by $m$ shows that this equals
$$ \sum_{h=1}^L  \sum_{m=l_1(h)}^{l_2(h) } (m^2  - 
 \epsilon(h)^2) r_{h,m} = \sum_{h=1}^L  \sum_{m=l_1(h)}^{l_2(h) } (\epsilon(h)-m)^2 r_{h,m} $$
$$ = \sum_{h=1}^L  \sum_{m=l_1(h)}^{l_2(h) } \sum_{i \in \Delta_{h,m}} \frac{(\epsilon(h)-m)^2}{|\!|e_i|\!|^2}
= \frac{1}{|\!|\beta |\!|^2},$$
and so 
$\sum_{(i,j) \in S} \lijb = 1.$
Finally note that
$$(\sum_{k_1 = l_1(h)}^m (\epsilon(h) - k_1) r_{h,k_1})(\sum_{k_2= l_1(h)}^{l_2(h)} r_{h,k_2})
= \sum_{k_1 = l_1(h)}^m \sum_{k_2= l_1(h)}^{l_2(h)} (k_2 r_{h,k_2} - k_1 r_{h,k_2})r_{h,k_1}.$$
This sum is positive because if $k_2 >m$ then the contribution of the pair $(k_1,k_2) $ to the
sum is $(k_2 - k_1) r_{h,k_1} r_{h,k_2} > 0$, whereas if $k_2 \leq m$ then the total contribution of
the pairs $(k_1,k_2)$ and $(k_2,k_1)$ is zero. Thus by (\ref{lam}) we have
$\lijb \geq 0$ for all $(i,j) \in S$, and hence $\beta $ lies in the convex hull of
$\{ e_i - e_j:(i,j) \in S\}$ as required.

\begin{lem} \label{5.3} If $\beta$ satisfies the conditions of Proposition
\ref{5.1} and if $i \in \Delta_{h,m}$ and $j \in \Delta_{h',m'}$, then
$$\beta.(e_i - e_j) = |\!| \beta |\!|^2$$
if and only if $h'=h$ and $m' = m+1$, and
$$\beta.(e_i - e_j) \geq |\!| \beta |\!|^2$$
if and only if either $m' \geq m+2$ or $m' = m+1$ and $h' \geq h$.
\end{lem}

\noindent{\bf Proof}: This follows immediately from the formula for 
$\beta$ and conditions on the function $\epsilon$ in the statement
of Proposition \ref{5.1}.

\begin{rem} \label{5.4}
Let $\beta$ be as in Proposition \ref{5.1}. Then there is a unique bijection
$\phi:J \to \{1,...,t\}$ from the indexing set $J$ of the partition
$\{ \Delta_{h,m}: (h,m) \in J \}$ of $\{1,...,M\}$ to the set of positive
integers $\{1,...,t\}$, where $t$ is the size of $J$, which takes the Hebrew
lexicographic ordering on $J$ to the standard ordering on integers; that is,
$\phi(h,m) \leq \phi(h',m')$ if and only if either $m<m'$ or $m=m'$ and
$h \leq h'$. We can define an increasing function
$$\delta: \{1,...,t\} \to \{1,...,t \}$$
such that $\delta(\phi(h,m))$ is the number of elements $(h',m') \in J$
such that either $m' < m+1$ or $m' = m+1$ and $h'<h$. Then
$\delta(k) \geq k$ for all $k \in \{1,...,t\}$, and if $(h,m)$ and $(h,m+1)$
both belong to $J$ then $\delta(\phi(h,m)) = \phi(h,m+1) -1$ and
$\delta(\phi(h,m)) < \delta(\phi(h,m)+1)$. Conversely if
$$k_2 -1 = \delta(k_1) < \delta(k_1 +1)$$
then there exists $(h,m) \in J$ with $(h,m+1) \in J$ such that
$k_1 = \phi(h,m)$ and $k_2 = \phi(h,m+1)$.

When it is helpful to make the dependence on $\beta$ explicit,
we shall write $\delta_\beta:\{1,...,t_\beta\} \to \{1,...,t_\beta\}$
and $\{ \delta_{h,m}(\beta):(h,m) \in J_\beta \}$.

Lemma \ref{5.3} tells us that if $i \in \Delta_{h,m}$ and $j \in \Delta_{h',m'}$
then
$\beta.(e_i - e_j) \geq |\!| \beta |\!|^2$
if and only if $\phi(h',m') > \delta(\phi(h,m))$.
\end{rem}

\begin{df} \label{5.5} Recall that if $1 \leq i \leq q$ then 
$e_{m_1 +...+m_{i-1} +1},...,e_{m_1 + ... + m_i}$ are the weights of the standard
representation on $\CC^{m_i}$ of the component $GL(m_i;\CC)$ of
$R_l = \prod_{i=1}^q GL(m_i;\CC)$. If $\beta$ and $\phi:J \to \{1,...,t\} $
are as in Remark \ref{5.4} and $1 \leq i \leq q$ and $1 \leq k \leq t$, then
set
$$\Delta^k = \Delta_{\phi^{-1}(k)}, \quad \Delta^k_i = \Delta_{\phi^{-1}(k)} \cap \{m_1 + ... + m_{i-1} + 1, ..., m_1 + ... + m_i \}$$
and let $m_i^k$ denote the size of $\Delta^k_i$, so that $m^1_i + ... + m_i^t = m_i$.
\end{df}

\begin{rem} \label{5.5a}
By Remark \ref{9.2} the partition $\{ \Delta^k(\beta) : 1 \leq k \leq t_\beta \}$ of
$\{1,...,M\}$ and the function $\delta_\beta: \{1,...,t_\beta\} \to \{1,...,t_\beta\}$ are determined
by $\beta$. Conversely, from the partition $\{ \Delta^k(\beta) : 1 \leq k \leq t_\beta \}$ of
$\{1,...,M\}$ and the function $\delta_\beta: \{1,...,t_\beta\} \to \{1,...,t_\beta\}$ we can recover
$\beta$ as the closest point to 0 of the convex hull of
$$\{e_i - e_j: i \in \Delta^{k_1}(\beta) \mbox{ and } j \in \Delta^{k_2}(\beta) \mbox{ where }
k_2 > \delta_\beta (k_1) \}.$$
Note, however, that although given any partition 
$\{ \Delta^k : 1 \leq k \leq t \}$ of
$\{1,...,M\}$ and increasing function $\delta: \{1,...,t\} \to \{1,...,t\}$ satisfying $\delta(k) \geq k$
for $1 \leq k \leq t$, we can consider  the closest point $\beta$ to 0 of the convex hull of
$$\{e_i - e_j: i \in \Delta^{k_1} \mbox{ and } j \in \Delta^{k_2} \mbox{ where }
k_2 > \delta (k_1) \},$$
it is not necessarily the case that the associated  partition $\{ \Delta^k(\beta) : 1 \leq k \leq t_\beta \}$ of
$\{1,...,M\}$ and function $\delta_\beta: \{1,...,t_\beta\} \to \{1,...,t_\beta\}$ coincide with the
given partition 
$\{ \Delta^k : 1 \leq k \leq t \}$ of
$\{1,...,M\}$ and function $\delta: \{1,...,t\} \to \{1,...,t\}$. For example, some amalgamation
and rearrangement may be needed as in the proof of Proposition 9.1.
\end{rem}

\section{Balanced $\delta$-filtrations}
\renorm

The last section studied the indexing set $\Gamma$ for the stratification
$\{ \Sigma_\gamma: \gamma \in \Gamma \}$ of $\C^{ss}$ defined as in $\S$2. 
In this section we will consider what it means for a semistable
holomorphic bundle over the Riemann surface $\Sigma$ to belong to
a stratum $\Sigma_\gamma = \Sigma_{\beta,l}$, where $\beta$ is as in Proposition \ref{5.1}.

\begin{df} \label{def5.8} 
We shall say that a semistable bundle $E$ has a $\delta$-filtration
$$0 = E_0 \subset E_1 \subset ... \subset E_t = E$$
with associated function $\delta:\{1,...,t\} \to \{1,...,t\}$ if $\delta$
is an increasing function such that  if $1 \leq k \leq t$ then
$\delta(k) \geq k$ and the induced filtration
$$0 \subset \frac{E_k}{E_{k-1}} \subset \frac{E_{k+1}}{E_{k-1}} \subset ... \subset 
\frac{E_{\delta(k)}}{E_{k-1}}$$
is trivial.
\end{df}

Let $G(\delta)$ be the graph with vertices $1,...,t$ and edges joining $i$ to $j$
if $j-1 = \delta(i) < \delta(i+1)$. Then the connected components of $G(\delta)$
are of the form 
$$ \{ i^h_{l_1(h)},...,i^h_{l_2(h)} \}$$
for $1 \leq h \leq L$, where $i^h_j - 1 = \delta(i^h_{j-1}) < \delta(i^h_{j-1} +1) $
if $1 < j \leq s_h$, and $i_1^h - 1 $ is not in the image of $\delta$ 
and either $\delta(i^h_{s_h}) = u$ or $\delta(i^h_{s_h}) =\delta(i^h_{s_h}+1) $.
Moreover $l_1(h) \leq l_2(h)$ can be chosen so that if
$$\epsilon(h) = \left( \sum_{m=l_1(h)}^{l_2(h)} m \tilde{r}_{h,m} \right)
\left( \sum_{m=l_1(h)}^{l_2(h)} \tilde{r}_{h,m} \right)^{-1},$$
where
$ \tilde{r}_{h,m} = \mbox{rank}(E_{i^h_m} / E_{i^h_m - 1}),$
then $-1/2 \leq \epsilon(h) < 1/2$, and the ordering of the
components of $G(\delta)$ can be chosen so that
\begin{equation} \label{eps}
 \epsilon(1) \geq \epsilon(2) \geq ... \geq \epsilon(L). \end{equation}
We shall say that the $\delta$-filtration is {\em balanced}
if the inequalities in (\ref{eps}) are all strict and if
\begin{equation} \label{Hebrew} i_{m_1}^{h_1} \leq i_{m_2}^{h_2} \quad
\mbox{ if and only if } \quad m_1 < m_2 \mbox{ or } m_1=m_2 \mbox{ and } h_1 \leq h_2;
\end{equation}
that is, if the usual ordering on $\{1,...,t\}$ is the same as the Hebrew
lexigraphic ordering via the pairs $(h,m)$.

\begin{rem} \label{sizeb}
If  $\beta$ is as in Proposition \ref{5.1} then the proof of that proposition
shows that
\begin{equation} \label{normb1}  \frac{1}{|\!|\beta|\!|^2} =
\sum_{h=1}^L \left( \frac{\sum_{l_1(h) \leq i < j \leq l_2(h)} (j-i)^2 r_{h,i} 
r_{h,j}}{\sum_{l_1(h) \leq i \leq l_2(h)} r_{h,i}} \right)
= \sum_{h=1}^L \sum_{m=l_1(h)}^{l_2(h)} (m - \epsilon (h))^2 r_{h,m} \end{equation}
where
$$r_{h,m} = \sum_{j \in \Delta_{h,m}} |\!| e_j|\!|^{-2} = \sum_{i=1}^q \sum_{j \in \Delta^{\phi(h,m)}_i}
|\!|e_j|\!|^{-2}$$
$$=  \sum_{i=1}^q \sum_{j \in \Delta^{\phi(h,m)}_i}
p_i =  \sum_{i=1}^q m^{\phi(h,m)}_i p_i
=  (1 -g+d/n) \sum_{i=1}^q m^{\phi(h,m)}_i n_i .$$
Thus 
\begin{equation} \label{normb2}  \frac{1}{|\!|\beta|\!|^2} = (1 -g + d/n)
\sum_{h=1}^L \left( \frac{\sum_{l_1(h) \leq i < j \leq l_2(h)} (j-i)^2 \tilde{r}_{h,i} 
\tilde{r}_{h,j}}{\sum_{l_1(h) \leq i \leq l_2(h)} \tilde{r}_{h,i}} \right)
\end{equation}
$$= (1 -g + d/n)  \sum_{h=1}^L \sum_{m=l_1(h)}^{l_2(h)} (m - \epsilon (h))^2 \tilde{r}_{h,m}$$
where
$$\tilde{r}_{h,m} = \sum_{i=1}^q m^{\phi(h,m)}_i n_i $$
and
$\epsilon(h)$ is given by
$$ \left( \sum_{m=l_1(h)}^{l_2(h)} m {r}_{h,m} \right)
\left( \sum_{m=l_1(h)}^{l_2(h)} {r}_{h,m} \right)^{-1}= \left( \sum_{m=l_1(h)}^{l_2(h)} m \tilde{r}_{h,m} \right)
\left( \sum_{m=l_1(h)}^{l_2(h)} \tilde{r}_{h,m} \right)^{-1}.$$
Note that if 
$0 = E_0 \subset E_1 \subset ... \subset E_t = E$
is a filtration such that 
$$E_k/E_{k-1} \cong \bigoplus_{i=1}^q \CC^{m_i^k} \otimes D_i$$
if $1 \leq k \leq t$,
where $t$ and $m_i^k$ for $1 \leq i \leq q$ and $1 \leq k \leq t$ are as in Definition
\ref{5.5} and $D_1,...,D_q$ are nonisomorphic stable bundles of ranks $n_1,...,n_q$
and all of the same slope $d/n$, then
$$\tilde{r}_{h,m} = \mbox{rank}(E_{\phi(h,m)}/ E_{\phi(h,m) -1}).$$
\end{rem}

\begin{prop} \label{5.6}
Let $\beta$ be as in Proposition \ref{5.1}, let $\delta$ be as in Remark \ref{5.4} 
and let $E$ be a semistable holomorphic
structure on $\mathcal{E}$.

(i) If $E$ represents an element of the stratum $\Sigma_{\beta,l}$ then $E$ has a 
unique balanced $\delta$-filtration
$0 = E_0 \subset E_1 \subset ... \subset E_t = E$
such that 
$$E_k/E_{k-1} \cong \bigoplus_{i=1}^q \CC^{m_i^k} \otimes D_i$$
and hence
$$E_{\delta(k)}/E_{k-1} \cong \bigoplus_{i=1}^q \CC^{m_i^k + ... + m^{\delta(k)}_i}
\otimes D_i$$
if $1 \leq k \leq t$,
where $t$ and $m_i^k$ for $1 \leq i \leq q$ and $1 \leq k \leq t$ are as in Definition
\ref{5.5} and $D_1,...,D_q$ are nonisomorphic stable bundles of ranks $n_1,...,n_q$
and all of the same slope $d/n$. 

(ii) Conversely, if $E$ has a 
balanced $\delta$-filtration
$0 = E_0 \subset E_1 \subset ... \subset E_t = E$
as in (i) then $E$ represents an element of the
stratum $\Sigma_{\beta,l}$ if and only if $E$ has no filtration with the corresponding
properties for any $\beta'$ satisfying $|\!| \beta' |\!| > |\!| \beta |\!|$.
\end{prop}

\noindent{\bf Proof}: Recall from (\ref{sgqy}) that
\begin{equation} \label{first} \Sigma_{\beta,l} =
 \G_c Y_{\beta,l}^{\backslash E}
\cong
 \G_c \times_{Q_\beta,l} Y_{\beta,l}^{\backslash E}
\end{equation}
and that
if $E$ represents an element of $Y_{\beta,l}^{\backslash E}$ then its
orbit under the complex one-parameter subgroup of $R_l$ generated
by $\beta$ has a limit point in $Z_{R_l}^s$. This limit point is represented by the bundle 
$\mbox{gr}(E)$ which is of
the form
$$ \mbox{gr}(E) \cong \bigoplus_{i=1}^q \CC^{m_i} \otimes D_i$$
where $D_1,...,D_q$ are nonisomorphic stable bundles of ranks $n_1,...,n_q$
and all of the same slope $d/n$. 
Recall also from \cite[$\S$ 7]{AB} that $\mathcal{C}$ is an infinite
dimensional affine space, and if we fix a $C^{\infty}$ identification of the
fixed $C^{\infty}$ hermitian bundle $\mathcal{E}$ with
$\bigoplus_{i=1}^q \CC^{m_i} \otimes D_i$
then we can identify $\mathcal{C}$ with the infinite dimensional vector space
$$\Omega^{0,1}(\mbox{End}(\bigoplus_{i=1}^q \CC^{m_i} \otimes D_i))$$
in such a way that the zero element of $\Omega^{0,1}(\mbox{End}(\bigoplus_{i=1}^q \CC^{m_i} \otimes D_i))$
corresponds to the given holomorphic structure on $\bigoplus_{i=1}^q \CC^{m_i} \otimes D_i$. 
With respect to this identification, the action of $R_l = \prod_{i=1}^q GL(m_i;\CC)$ on
$\mathcal{C}$ is the action induced by the obvious action of $R_l$ on
$\bigoplus_{i=1}^q \CC^{m_i} \otimes D_i$. The one-parameter subgroup of $R_l$
generated by $\beta$ acts diagonally on $\CC^{m_1} \oplus ... \oplus \CC^{m_q}$ with
weights $\beta.e_j$ for $j \in \{1,...,M\}$ where $M = m_1 +...+m_q$, and so it acts on
$$\Omega^{0,1}(\mbox{End}(\bigoplus_{i=1}^q \CC^{m_i} \otimes D_i))
= \bigoplus_{i_1,i_2=1}^q \Omega^{0,1} (\CC^{m_{i_1}} \otimes (\CC^{m_{i_2}})^* \otimes D_{i_1}
\otimes D_{i_2}^*)$$
with weights $\beta.(e_i - e_j)$ for $i,j \in \{1,...,M\}$. If $E \in \Sigma_{\beta}$ then $E$ lies
in the $\mathcal{G}_c$-orbit of an element of the sum of those weight spaces for which the
weight $\beta.(e_i - e_j)$ satisfies $\beta. (e_i - e_j) \geq |\!| \beta |\!|^2$. By Remark \ref{5.4}
and Definition \ref{5.5} we have a partition $\{\Delta^1,...,\Delta^t\}$ of $\{1,...,M\}$ such
that  $\beta. (e_i - e_j) \geq |\!| \beta |\!|^2$ if and only if $i \in \Delta^{k_1}$ and
$j \in \Delta^{k_2}$ where $k_2 > \delta(k_1)$. So if we make identifications
$$ \CC^M = \bigoplus_{i=1}^q \CC^{m_i} = \bigoplus_{i=1}^q \bigoplus_{k=1}^t \CC^{m_i^k}$$
using the induced partition $\{ \Delta^k_i: 1 \leq i \leq q, 1 \leq k \leq t \}$ of 
$\{ 1,...,M\}$ as in Definition \ref{5.5}, then any $E \in \Sigma_\beta$ lies in the
$\mathcal{G}_c$-orbit of an element of
$$\bigoplus_{i_1,i_2=1}^q \bigoplus_{k_1=1}^t \bigoplus_{k_2 = \delta(k_1) + 1}^t
\Omega^{0,1} (\CC^{m^{k_1}_{i_1}} \otimes (\CC^{m^{k_2}_{i_2}})^* \otimes D_{i_1}
\otimes D_{i_2}^*).$$
This completes the proof of (i), as such an element of
$\Omega^{0,1}(\mbox{End}(\bigoplus_{i=1}^q \CC^{m_i} \otimes D_i))$
represents a holomorphic structure $E$ on $\mathcal{E}$ with a filtration of the
required form (uniqueness follows from (\ref{first}) and the fact that
$Q_{\beta,l}$ preserves the filtration), and (ii) is now a consequence of (\ref{sb}).

\begin{cor} \label{nice}
Let $\beta$ be as in Proposition \ref{5.1}, let $\delta$ be as in Remark \ref{5.4} 
and let $E$ be a semistable holomorphic
structure on $\mathcal{E}$ with a 
balanced $\delta$-filtration
\begin{equation} \label{niceq} 0 = E_0 \subset E_1 \subset ... \subset E_t = E \end{equation}
whose subquotients satisfy
$$E_k/E_{k-1} \cong \bigoplus_{i=1}^q \CC^{m_i^k} \otimes D_i,$$
where $t$ and $m_i^k$ for $1 \leq i \leq q$ and $1 \leq k \leq t$ are as in Definition
\ref{5.5} and $D_1,...,D_q$ are nonisomorphic stable bundles of ranks $n_1,...,n_q$
and all of the same slope $d/n$.  Then $E$ represents an element of the stratum
$\Sigma_{\beta,l}$ if and only if, in the notation of Definition \ref{def5.8}, there is no
$h \in \{1,...,L \}$ and refinement
$$0 = E_0 \subset ... \subset E_{i^h_{l_1(h)}-1} \subset F_{l_1(h)} \subset
E_{i^h_{l_1(h)}} \subset ... 
\subset E_{i^h_m-1} \subset F_m \subset
E_{i^h_m} \subset ... $$
$$... \subset E_{i^h_{l_2(h)}-1} \subset F_{l_2(h)} \subset
E_{i^h_{l_2(h)}} \subset ... 
\subset E_t = E$$
of (\ref{niceq}) with
$$\frac{\sum_{m=l_1(h)}^{l_2(h)} m \mbox{ {\rm rank}}(F_m/E_{i^h_m -1})}
{\sum_{m=l_1(h)}^{l_2(h)}  \mbox{ {\rm rank}}(F_m/E_{i^h_m -1})} <
\frac{\sum_{m=l_1(h)}^{l_2(h)} m \mbox{ {\rm  rank}}(E_{i^h_m}/E_{i^h_m -1})}
{\sum_{m=l_1(h)}^{l_2(h)}  \mbox{ {\rm rank}}(E_{i^h_m}/E_{i^h_m -1})} 
=\epsilon(h),$$
such that the induced filtrations
$$0 \subset \frac{E_{i^h_m -1}}{E_{i^{h_1}_{m_1}-1}} \subseteq
\frac{F_m}{E_{i^{h_1}_{m_1}-1}}$$
with 
$$m_1 - \epsilon(h_1) \leq \frac{\sum_{m=l_1(h)}^{l_2(h)} m 
\mbox{ {\rm rank}}(F_m/E_{i^h_m -1})}
{\sum_{m=l_1(h)}^{l_2(h)}  \mbox{ {\rm rank}}(F_m/E_{i^h_m -1})} $$
and
$$ 0\subset \frac{E_{i^h_m}}{F_m} \subseteq \frac{E_{i^{h_2}_{m_2}}}{F_m}$$
with
$$m_2 - \epsilon(h_2) \geq m - \frac{\sum_{m=l_1(h)}^{l_2(h)} m 
\mbox{ {\rm rank}}(F_m/E_{i^h_m -1})}
{\sum_{m=l_1(h)}^{l_2(h)}  \mbox{ {\rm rank}}(F_m/E_{i^h_m -1})} $$
are all trivial.
\end{cor}

\begin{rem}
If $ 0 = E_0 \subset E_1 \subset ... \subset E_{j-1} \subset F \subset E_j ... \subset E_t = E$
is a refinement of the filtration (\ref{niceq}) of $E$ such that the induced filtration
$$ 0 \subseteq \frac{E_{j-1}}{E_i} \subset \frac{F}{E_i}$$
is trivial for some $i < j-1$, then $F/E_i$ is isomorphic to
$$ \frac{E_{j-1}}{E_i} \oplus \frac{F}{E_{j-1}}$$
and so $E$ can be given a filtration of the form
$$ 0 = E_0 \subset E_1 \subset ... \subset E_i \subset E'_{i+1} \subset E'_{i+2}
\subset ...
\subset E'_{j-1} \subset F \subset E_j \subset ... \subset E_t = E$$
where $E'_{i+1}/E_i \cong F/E_{j-1}$, $E'_{k+1}/E'_k \cong E_k/E_{k-1}$ for
$i<k<j$ and $F/E'_{j-1} \cong E_{j-1}/E_{j-2}$. A similar result is true
if the induced filtration
$$0 \subset \frac{F}{E_{j-1}} \subset \frac{E_i}{E_{j-1}}$$
is trivial for some $i>j$.
\end{rem}

\noindent{\bf Proof of Corollary \ref{nice}}: This follows from \cite{K7} and
the proof of Proposition \ref{5.1}, which tells us that if $\beta'\neq \beta$ is the closest 
point to 0 of the convex hull of $\{e_i - e_j:(i,j) \in S' \}$ where $S'$ is a 
subset of $S$, then $S'$ can be chosen so that the connected components
of the graph $G(S')$ give us a refinement of the partition $\{\Delta_{h,m}:(h,m) \in J\}$
of $\{1,...,M\}$ associated to $\beta$, which in turn gives us a refinement of
the filtration (\ref{niceq}) with the required properties.

\begin{rem} 
Recall from Proposition \ref{5.6} that a semistable bundle $E$ represents an element of the
stratum $\Sigma_{\beta,l}$ if and only if if has a balanced $\delta$-filtration 
\begin{equation} \label{honk}
0 = E_0 \subset E_1 \subset ... \subset E_t = E
\end{equation}
such that if $1 \leq k \leq t$ then
$E_k/E_{k-1}$ is  of the form $ \bigoplus_{i=1}^q \CC^{m_i^k} \otimes D_i$
where $D_1,...,D_q$ are nonisomorphic stable bundles of ranks $n_1,...,n_q$
and all of the same slope $d/n$, and moreover $E$ has no balanced
$\delta$-filtration with the corresponding properties
when $\beta$ is replaced with $\beta'$ satisfying $|\!| \beta' |\!| > |\!| \beta |\!|$.
From (\ref{normb2}) we have
\begin{equation} \label{normb}  \frac{1}{|\!|\beta|\!|^2} = (1-g+d/n)
\sum_{h=1}^L \left( \frac{\sum_{l_1(h) \leq i < j \leq l_2(h)} (j-i)^2 \tilde{r}_{h,i} 
\tilde{r}_{h,j}}{\sum_{l_1(h) \leq i \leq l_2(h)} \tilde{r}_{h,i}} \right),
\end{equation}
where 
$$\tilde{r}_{h,m} = \mbox{rank}(E_{\phi(h,m)}/ E_{\phi(h,m) -1}).$$
This gives us some sort of measure of the triviality of the 
balanced $\delta$-filtration
(\ref{honk}); very roughly speaking, the more trivial this filtration, the
smaller the size of $|\!|\beta|\!|^{-2}$ and hence the larger $|\!| \beta|\!|$ becomes.

Let us therefore define the {\em triviality} of the balanced
$\delta$-filtration (\ref{honk}) with 
associated function $\delta$ to be
$$ \left( \sum_{h=1}^L \left( \frac{\sum_{l_1(h) \leq i < j \leq l_2(h)} (j-i)^2 \tilde{r}_{h,i} 
\tilde{r}_{h,j}}{\sum_{l_1(h) \leq i \leq l_2(h)} \tilde{r}_{h,i}} \right) \right)^{-1/2}$$
\begin{equation} \label{triv}
=\left( \sum_{h=1}^L \sum_{m=l_1(h)}^{l_2(h)} (m - \epsilon(h))^2 \tilde{r}_{h,m} \right)^{-1/2};
\end{equation}
Remarks \ref{9.2} and \ref{5.5a} tell us that this is well defined.
Thus if $E \in \Sigma_{\beta,l}$ the balanced $\delta$-filtration (\ref{honk}) 
associated to $E$ by Proposition 10.1 can be thought of
as having maximal triviality (according to this measure) among the
balanced $\delta$-filtrations of $E$.
\end{rem}

\begin{rem} \label{5.6a} 
Let $\beta$ be as in Proposition \ref{5.1}, let $\delta$ be as in Remark 
\ref{5.4}  and let $E$ be a semistable holomorphic
structure on $\mathcal{E}$ with a $\delta$-filtration 
$$0 = E_0 \subset E_1 \subset ... \subset E_t = E$$
such that if $1 \leq k \leq t$ then
$E_k/E_{k-1} \cong \bigoplus_{i=1}^q \CC^{m_i^k} \otimes D_i$
where $D_1,...,D_q$ are nonisomorphic stable bundles of ranks $n_1,...,n_q$
and all of the same slope $d/n$. Then the proof of Proposition \ref{5.6} 
and (\ref{stratclos}) shows that $E$ represents an element of $\Sigma_{\tilde{\gamma}}$
for some $\tilde{\gamma} \in \Gamma$ satisfying $\tilde{\gamma} \geq \gamma = (\beta,l)$ with respect to the
partial ordering on $\Gamma$ defined just before (\ref{stratclos}).
\end{rem}

\begin{rem} \label{5.6b}
If $[{\bf n,m}] \in \mathcal{I}^{ss}$ then we can define $\beta [{\bf n,m}] \in \Gamma$ as
follows. Let $m_i = \sum_{j=1}^r m_{ij}$ and suppose that $R_l$ is the
subgroup of $GL(p;\CC)$ given by the embedding of $\prod_{i=1}^q GL(m_i;\CC)$
via a fixed identification of $\bigoplus_{i=1}^q \CC^{n_i(1-g+d/n)} \otimes
\CC^{m_i}$ with $\CC^p$. If $q=1$ let 
$$\beta [{\bf n,m}] = R_l,$$
and if $q>1$ let $\beta [{\bf n,m}] \in \mathcal{B}_l \backslash \{ 0 \}$ be
represented by the closest point to 0 of the convex hull of the weights of the
representation of the maximal compact torus $T_l$ of $R_l$ on
$$\bigoplus_{i_1, i_2 = 1}^q \bigoplus_{1 \leq j_1 < j_2 \leq r} \CC^{m_{i_1 j_1}}
\otimes (\CC^{m_{i_2 j_2}})^*$$
given by identifying $\bigoplus_{j=1}^r \CC^{m_{ij}}$ with $\CC^{m_i}$ for
$1 \leq i \leq q$. Then we have from Remark \ref{5.6a} that
\begin{equation} \label{compare1}
S_{[{\bf n,m}]} \subseteq \bigcup_{\gamma \geq \beta [{\bf n,m}]} \Sigma_\gamma,
\end{equation}
and from Remark \ref{4.3a} and Proposition \ref{5.6} that
\begin{equation} \label{compare2}
\Sigma_{\beta [{\bf n,m}]} \subseteq \bigcup_{[{\bf n',m'}] \geq [{\bf n,m}]}
S_{[{\bf n',m'}]}
\end{equation}
where $\geq$ denotes in (\ref{compare1}) the partial order on $\Gamma$ used in
Remark \ref{5.6a}, whereas in (\ref{compare2}) it denotes the partial order on
$\mathcal{I}^{ss}$ described in Remark \ref{4.3a}.
\end{rem}

\begin{rem} \label{5.7} It follows from Proposition \ref{5.6} that if $R_l$ is as at
(\ref{rl}) then a holomorphic structure belongs to
\[
\bigcup_{\beta \in \mathcal{B}_l \backslash \{0\}} \Sigma_{\beta,l}
\]
if and only if 
$ E \not \cong \mbox{gr}(E) \cong (\CC^{m_1} \otimes D_1) 
\oplus \cdots \oplus (\CC^{m_q} \otimes D_q) $
where $D_1, \ldots ,D_q$ are all stable of slope $d/n$ and ranks $n_1, \ldots, n_q$ and are not isomorphic 
to one another.
\end{rem}

\section{Pivotal filtrations}
\renorm

Let us now consider the relationship between the balanced $\delta$-filtration
associated to a semistable bundle $E$ as in Proposition \ref{5.6} and the
maximal and minimal Jordan--H\"{o}lder filtrations defined in $\S$8.

Indeed, motivated by Proposition \ref{5.6}, we can try to carry our analysis of the maximal
Jordan--H\"{o}lder filtration
\begin{equation} \label{habove} 0=E_0 \subset E_1 \subset ... \subset E_t = E \end{equation}
of a bundle $E$ a bit further. Recall that if $1 \leq j \leq t$ then the 
subquotient $E_j/E_{j-1}$ is the maximal subbundle of $E/E_{j-1}$ which is
a direct sum of stable bundles all having maximal slope among the nonzero
subbundles of $E/E_{j-1}$. We can ask whether it is true for every subbundle $F$ 
of $E$ satisfying $E_{j-1} \subset F \subset E_j$ and $\mbox{slope}(E_j /F)
= \mbox{slope}(E_j/E_{j-1})$ that $E_j/F$ is the maximal subbundle of
$E/F$  which is
a direct sum of stable bundles all having maximal slope among the nonzero
subbundles of $E/F$. Of course if $E_j/E_{j-1}$ is itself stable there are no
such intermediate subbundles $F$, so this is trivially true, but it is not
always the case (as Example \ref{5.8} below shows). 

If there does
exist such an intermediate subbundle $F$,
then by Lemma \ref{e4.2} 
both $F/E_{j-1}$ and $E_j/F$ are of the form
\begin{equation} \label{above}
(\CC^{m_1} \otimes D_1) \oplus \cdots \oplus (\CC^{m_q} \otimes D_q)
\end{equation}
where $D_1, \ldots ,D_q$ are all stable of slope $d/n$ and ranks $n_1, \ldots, n_q$ and are not isomorphic 
to one another. So
we can then ask whether it is possible to find a canonical refinement 
\begin{equation} \label{tabove} 0 = F_0 \subset F_1 \subset ... \subset F_u = E \end{equation}
of the maximal Jordan--H\"{o}lder filtration (\ref{habove}) of $E$
with an increasing function $\delta:\{1,...,u\} \to \{1,...,u\}$ 
such that $\delta(j) \geq j$ if $1 \leq j \leq u$, each subquotient
$F_{\delta(j)}/F_{j-1}$ is of the form (\ref{above}) and
moreover for every subbundle $F$ of $E$ satisfying $F_{j-1} \subseteq F \subset F_j$
and $\mbox{slope}(F_j/F) = \mbox{slope}(F/F_{j-1})$, the quotient
$F_{\delta(j)}/F$ is the maximal subbundle of $E/F$  which is 
a direct sum of stable bundles all having maximal slope among the nonzero
subbundles of $E/F$. 
The following example shows that even this cannot be achieved in a canonical way.

\begin{example} \label{5.8} Let $E_1$ and $E_2$ be semistable bundles
over $\Sigma$ such that 
$\mbox{slope}(E_1)$ equals $\mbox{slope}(E_2)$, with
maximal Jordan--H\"{o}lder filtrations
$$0 \subset D_1 \subset E_1 \quad \mbox{ and } \quad 0 \subset D_2 \subset E_2$$
where $D_1$, $D_2$, $E_1/D_1$ and $E_2/D_2$  are nonisomorphic 
stable bundles all of the same slope as $E_1$ and $E_2$, and let
$$E = E_1 \oplus E_2.$$
We observed in Remark \ref{4.7} that the maximal Jordan--H\"{o}lder filtration
of a direct sum of semistable bundles of the same slope is the direct
sum of their maximal Jordan--H\"{o}lder filtrations (with the shorter one extended
trivially at the top if they are not of the same length). Thus the maximal
Jordan--H\"{o}lder filtration of $E$ is
\begin{equation} \label{filt1} 0 \subset D_1 \oplus D_2 \subset E_1 \oplus E_2 = E.
\end{equation}
By Lemma \ref{e4.2} there are precisely two proper subbundles $F$ of
$D_1 \oplus D_2$ with $\mbox{slope}(F) = \mbox{slope}(D_1 \oplus D_2 / F) =
\mbox{slope}(D_1 \oplus D_2) $, namely $D_1$ and $D_2$. The
maximal Jordan--H\"{o}lder filtration of $E/D_1 = (E_1/D_1) \oplus E_2$ is
$$ 0 \subset (E_1/D_1) \oplus D_2 \subset (E_1/D_1) \oplus E_2,$$
so we can refine the filtration (\ref{filt2}) of $E$ to get
\begin{equation} \label{filt2} 0 = F_0 \subset F_1 \subset F_2 \subset F_3 \subset F_4 = E
\end{equation}
where $F_1 = D_1$, $F_2 = D_1 \oplus D_2$ and $F_3 = E_1 \oplus D_2$, and we can 
define $\delta : \{1,2,3,4\} \to \{1,2,3,4\}$ by $\delta(1) = 2$, $\delta(2) = 3$
and $\delta(3) = \delta(4) = 4$. If $1\leq j \leq 4$ then $F_{\delta(j)}/F_{j-1}$ is the maximal subbundle
of $E/F_{j-1}$ which is a direct sum of stable subbundles all having maximal slope
among the nonzero subbundles of $E/F_{j-1}$. Moreover this is trivially still true if we 
replace $F_{j-1}$ by any subbundle $F$ of $E$ satisfying $F_{j-1} \subseteq F \subset F_j$
and $\mbox{slope}(F_j/F) = \mbox{slope}(F_j/F_{j-1})$, since the only such $F$ is $F_{j-1}$
itself.
We can of course reverse the r\^{o}les of $E_1$ and $E_2$ in this construction, to
get another refinement 
\begin{equation} \label{filt3} 0  \subset  D_2 \subset D_1 \oplus D_2 \subset D_1\oplus E_2 \subset 
E_1 \oplus E_2 = E
\end{equation}
of (\ref{filt1}). Thus there are precisely two refinements of the maximal Jordan--H\"{o}lder
filtration of $E_1 \oplus E_2$ with the required propertes, and if $E_1$ has
the same rank as $E_2$  and $D_1$ has the same rank as $D_2$ then by symmetry
there can be no canonical choice.

Notice that if $\mbox{rank}(D_1)/\mbox{rank}(E_1) = \mbox{rank}(D_2)/\mbox{rank}(E_2)$
then neither of the $\delta$-filtrations (\ref{filt2}) and (\ref{filt3}) is balanced since the
inequalities (\ref{eps}) are not strict; 
however the $\delta$-filtration (\ref{filt1}) is balanced and has maximal triviality,
in the sense of (\ref{triv}), among balanced $\delta$-filtrations of $E_1 \oplus E_2$.
If on the other hand 
$\mbox{rank}(D_1)/\mbox{rank}(E_1) \neq \mbox{rank}(D_2)/\mbox{rank}(E_2)$
then precisely one of the $\delta$-filtrations (\ref{filt2}) and (\ref{filt3}) is balanced
and it has maximal triviality, in the sense of (\ref{triv}), among balanced $\delta$-filtrations
of $E_1 \oplus E_2$. This filtration then determines the stratum $\Sigma_\gamma$ to which
$E$ belongs, and in this case $E$ represents an element of the open subset
$\Sigma^s_{\beta,l}$  
of $\Sigma_\gamma$. 

\end{example}


\begin{lem} \label{5.9} Let $E$ be a bundle over $\Sigma$ with a filtration
$$0 = F_0 \subset F_1 \subset ... \subset F_u = E$$
and let $\{\Delta_1,...,\Delta_L\}$ be a partition of $\{1,...,u\}$ such that if $1 \leq h \leq L$
and $\Delta_h = \{i^h_1,...,i^h_{s_h} \}$ where $i^h_1 < i^h_2 < ... < i^h_{s_h}$, then the induced extension
\begin{equation} \label{ijh} 
0 \to  F_{i^h_j}/F_{i^h_j-1} \to  F_{i^h_{j+1} - 1}/F_{i^h_j-1} \to F_{i^h_{j+1} - 1}/F_{i^h_j} \to 0
\end{equation}
is trivial. Then we can associate to this filtration of $E$, partition $\{\Delta_1,...,\Delta_L\}$ of $\{1,...,u\}$
and trivialisations of the induced extensions (\ref{ijh}) a sequence of elements of
$$H^1( \Sigma,  (F_{i^h_j}/F_{i^h_j-1}) \otimes (F_{i^h_{j+1}}/F_{i^h_{j+1} -1})^*   )$$
or equivalently of extensions
$$
0 \to  F_{i^h_j}/F_{i^h_j-1} \to  E^h_j \to F_{i^h_{j+1}}/F_{i^h_{j+1} -1} \to 0
$$
for $1 \leq h \leq L$ and $1 \leq j \leq s_h -1$.
\end{lem}

\noindent{\bf Proof}: This lemma follows immediately from the well known bijective
correspondence between holomorphic extensions of a holomorphic bundle $D_1$ over
$\Sigma$ by another holomorphic bundle $D_2$ and elements of $H^1(\Sigma, D_1 \otimes
D_2)$. The extension
$$
0 \to  F_{i^h_{j+1} -1}/F_{i^h_j-1} \to  F_{i^h_{j+1}}/F_{i^h_j-1} \to F_{i^h_{j+1} }/F_{i^h_{j+1} -1} \to 0
$$
induced by the given filtration gives us an element of
$$H^1( \Sigma,  (F_{i^h_{j+1} -1}/F_{i^h_j-1}) \otimes (F_{i^h_{j+1}}/F_{i^h_{j+1} -1})^*   )$$
and the given trivialisation of the  extension (\ref{ijh}) gives us a decomposition of this as
$$H^1( \Sigma,  (F_{i^h_j}/F_{i^h_j-1}) \otimes (F_{i^h_{j+1}}/F_{i^h_{j+1} -1})^*   )
\oplus H^1( \Sigma,  (F_{i^h_{j+1} -1}/F_{i^h_j}) \otimes (F_{i^h_{j+1}}/F_{i^h_{j+1} -1})^*   ).$$
Projection onto the first summand gives us an extension
$$0 \to  F_{i^h_j}/F_{i^h_j-1} \to  E^h_j \to F_{i^h_{j+1}}/F_{i^h_{j+1} -1} \to 0
$$
as required.

\begin{rem} \label{5.10}
Lemma \ref{5.3} and Remark \ref{5.4} tell us that if $i \in \Delta_{i_1}^{k_1}$ and
$j\in \Delta_{i_2}^{k_2}$  where $k_2 > \delta(k_1)$ then $(e_i - e_j).  \beta
\geq |\!| \beta |\!|^2$, and equality occurs if and only if there exists $(h,m) \in J$ with
$(h,m+1) \in J$ such that $k_1 = \phi(h,m)$ and $k_2 = \phi(h,m+1)$. Thus,
retaining the notation of the proof of Proposition \ref{5.6}, we observe that if $E$ is represented
by an element of
$$\bigoplus_{i_1,i_2=1}^q \bigoplus_{k_1=1}^t \bigoplus_{k_2 = \delta(k_1) + 1}^t
\Omega^{0,1} (\CC^{m^{k_1}_{i_1}} \otimes (\CC^{m^{k_2}_{i_2}})^* \otimes D_{i_1}
\otimes D_{i_2}^*)$$
then the limit in $\mathcal{C}$ as $t \to \infty$ of $\exp(-it\beta) E$ is the bundle
$$ \mbox{gr} (E) \cong \bigoplus_{i=1}^q \CC^{m_i} \otimes D_i$$
which is represented by the zero vector in
$$\bigoplus_{i_1,i_2=1}^q \bigoplus_{k_1=1}^t \bigoplus_{k_2 = \delta(k_1) + 1}^t
\Omega^{0,1} (\CC^{m^{k_1}_{i_1}} \otimes (\CC^{m^{k_2}_{i_2}})^* \otimes D_{i_1}
\otimes D_{i_2}^*),$$
and so the limit $p_\beta (E)$ of $\exp(-it \beta )E$ in the blow-up of $\mathcal{C}$
along $\mathcal{G}_c Z_{R_l}^{ss}$ is an element of the fibre
$$\PP(\mathcal{N}_{l, \mbox{gr} E}) = \PP(
H^1(\Sigma, \bigoplus_{i_1,i_2=1}^{q} 
\CC^{m_{i_1}m_{i_2}-\delta_{i_1}^{ i_2}}  \otimes D_{i_1} \otimes D^*_{i_2})
$$
of the exceptional divisor over $\mbox{gr} E$. Indeed $p_\beta (E)$ is the element of
this fibre represented by the sum in
$$\bigoplus_{h=1}^L \bigoplus_{m=l_1(h)}^{l_2(h)-1}
H^1(\Sigma, (\bigoplus_{i_1}^{q} 
\CC^{m_{i_1}^{\phi(h,m)}} \otimes D_{i_1}) \otimes (\bigoplus_{i_2}^q 
\CC^{ m_{i_2}^{\phi(h,m+1)}} \otimes D_{i_2})^*)
$$
of the elements of 
$$
H^1(\Sigma, (\bigoplus_{i_1}^{q} 
\CC^{m_{i_1}^{\phi(h,m)}} \otimes D_{i_1}) \otimes (\bigoplus_{i_2}^q 
\CC^{ m_{i_2}^{\phi(h,m+1)}}  \otimes D_{i_2})^*)$$ for $1 \leq h \leq L$ and $l_1(h) \leq
m \leq l_2(h)$
corresponding to the extensions
$$0 \to \bigoplus_{i_1}^{q} 
\CC^{m_{i_1}^{\phi(h,m)}} \otimes D_{i_1} \to E^h_m \to
\bigoplus_{i_2}^q 
\CC^{ m_{i_2}^{\phi(h,m+1)}}  \otimes D_{i_2} \to 0$$
associated as in Lemma \ref{5.9} to the $\delta$-filtration
$ 0=E_0 \subset E_1 \subset ... \subset E_t = E$ of Proposition
\ref{5.6}. 
\end{rem}

\begin{prop} \label{5.11}
Let $\beta$ be as in Proposition \ref{5.1} and let $E$ be a semistable bundle
representing an element of $\Sigma_{\beta,l}$ with a balanced $\delta$-filtration
$$0=E_0 \subset E_1 \subset ... \subset E_t = E$$
such that if $1 \leq k \leq t$ then $E_k/E_{k-1} \cong \bigoplus_{i=1}^q \CC^{m_i^k} \otimes D_i$
as in Proposition \ref{5.6}. Suppose also that $k_1 \in \{1,...,t\}$ is
such that
$$\beta.e_j <0 \mbox{ whenever } j \in \Delta^{k_2} \mbox{ with } k_2 > \delta(k_1).$$
Then whenever $F$ is a subbundle of $E$ with $\mbox{slope}(F)= \mbox{slope}(E)$
and such that $E_{{k_1}-1} \subseteq F \subset E_{k_1}$, the subquotient $E_{\delta({k_1})}/F$
is the maximal subbundle of $E/F$ which is a direct sum of stable bundles all
having the same slope as $E/F$.
\end{prop}

\noindent{\bf Proof}: Since $F/E_{{k_1}-1}$ is a subbundle of 
$E_{k_1}/E_{{k_1}-1} \cong \bigoplus_{i=1}^q \CC^{m_i^{k_1}} \otimes D_i$
having the same slope as $D_1,...,D_q$, it follows from Lemma
\ref{e4.2} that
$$F/E_{{k_1}-1} = \bigoplus_{i=1}^q U_i \otimes D_i$$
where $U_i$ is a linear subspace of $\CC^{m_i^{k_1}}$ for $1 \leq i \leq q$, and so
$$E_{\delta({k_1})}/F \cong \bigoplus_{i=1}^q ((\CC^{m_i^{k_1}}/U_i) \oplus \CC^{m_i^{{k_1}+1} + ... +
m_i^{\delta({k_1})}}) \otimes D_i$$
is a sum of stable bundles all having the same slope (which is equal to $\mbox{slope}(E)$
and $\mbox{slope}(E/F)$). Let us suppose for a contradiction that $E/F$ has a subbundle
$E'/F$ which is not contained in $E_{\delta({k_1})}/F$ and which is of the required
form. Then we can choose $k_2>\delta(k_1)$ such that $E' \subseteq E_{k_2}$ but $E'$
is not contained in $E_{k_2 -1}$, and then the inclusion of $E'$ in $E_{k_2}$ induces a
nonzero map
$$ \theta: E'/F \to E_{k_2}/E_{k_2 -1} \cong \bigoplus_{i=1}^q \CC^{m_i^{k_2}} \otimes D_i.$$
Since nonzero bundle maps between stable bundles of the same slope are always isomorphisms,
by replacing $E'$ by a suitable subbundle we can assume that $E'/F \cong D_{i_0}$ for some
$i_0 \in \{1,...,q\}$, and that we can decompose $\CC^{m_{i_0}^{k_2}}$ as
$\CC\oplus \CC^{m_{i_0}^{k_2}-1} $ in such a way that the projection $\theta_0:E'/F \to D_{i_0}$
of $\theta$ onto the corresponding component $D_{i_0}$ of $E_{k_2}/E_{k_2 -1}$ is an
isomorphism. Then
$$\theta_0^{-1}:D_{i_0} \to E'/F \subseteq E_{k_2}/F$$
gives us a trivialisation of the extension of $E_{k_2-1}/E_{k_1}$ by this component 
$D_{i_0}$ of $E_{k_2}/E_{k_2 -1}$. By the definition of $\Sigma_{\beta,l}$ the limit $p_\beta(E) \in
\PP(\mathcal{N}_{l, \mbox{gr}E})$
of $\exp(-it\beta)E$ as $t \to \infty$ is semistable for the induced action of $\stab(\beta)/T^c_\beta$
where $\stab(\beta)$ is the stabiliser of $\beta$ under the coadjoint action
of $R_l$ and $T^c_\beta$ is the complex subtorus generated by $\beta$ (see \cite{K7}),
and by Remark \ref{5.10} $p_\beta(E)$ is represented by the sum over $h \in \{1,...,L\}$
and $m \in \{l_1(h),...,l_2(h)\}$ of the elements of 
$$
H^1(\Sigma, (\bigoplus_{i_1}^{q} 
\CC^{m_{i_1}^{\phi(h,m)}} \otimes D_{i_1}) \otimes (\bigoplus_{i_2}^q 
\CC^{ m_{i_2}^{\phi(h,m+1)}}  \otimes D_{i_2})^*)$$ 
corresponding to the extensions
$$0 \to \bigoplus_{i_1}^{q} 
\CC^{m_{i_1}^{\phi(h,m)}} \otimes D_{i_1} \to E^h_m \to
\bigoplus_{i_2}^q 
\CC^{ m_{i_2}^{\phi(h,m+1)}}  \otimes D_{i_2} \to 0$$
associated by Lemma \ref{5.9} to the filtration
$ 0=E_0 \subset E_1 \subset ... \subset E_t = E$.

Let $S_0$ be the set of ordered pairs $(i,j)$ with $i,j \in \{1,...,M\}$
such that the component of $p_\beta(E)$ in the weight space corresponding
to the weight $e_i - e_j$ for the action of the maximal torus $T_l$ of
$R_l = \prod_{i=1}^q GL(m_i; \CC)$ on
$$\bigoplus_{h=1}^L \bigoplus_{m=l_1(h)}^{l_2(h)-1}
H^1(\Sigma, (\bigoplus_{i_1}^{q} 
\CC^{m_{i_1}^{\phi(h,m)}} \otimes D_{i_1}) \otimes (\bigoplus_{i_2}^q 
\CC^{ m_{i_2}^{\phi(h,m+1)}} \otimes D_{i_2})^*)
$$
is nonzero. Since $p_\beta(E)$ is semistable for the action of $\stab(\beta)/T_\beta^c$,
it follows that $\beta$ is the closest point to 0 in the convex hull  of $\{e_i - e_j:(i,j) \in S_0\}$.
We may assume that $T_l$ acts diagonally with respect to the decomposition
of $\CC^{m_{i_0}^{k_2}}$ as $\CC \oplus \CC^{m_{i_0}^{k_2}-1}$; let $e_{j_0}$ be
the weight of the action of $T_l$ on the component $\CC$ of $\CC^{m_{i_0}^{k_2}}$ with
respect to this decomposition. Since $k_2 > \delta(k_1)$ and the extension of $E_{k_2-1}/E_{k_1}$
by the component of 
$E_{k_2}/E_{{k_2}-1} \cong \bigoplus_{i=1}^q \CC^{m_i^{k_2}} \otimes D_i$
corresponding to the weight $e_{j_0}$ is trivial, it follows that if $(i,j) \in S_0$ then
$j \neq j_0$. Since $\beta$ lies in the convex hull of $\{e_i - e_j:(i,j) \in S_0\}$
and $e_1,...,e_M$ are mutually orthogonal, this means that
$$\beta.e_{j_0} \geq 0,$$
and as $j_0 \in \Delta^{k_2}$ and $k_2 > \delta(k_1)$,  this gives us the required
contradiction.

\begin{rem} \label{5.12} Dual to the definition of $\delta$ in Remark \ref{5.4},
we can define an increasing function $\delta': \{1,...,t\} \to \{1,...,t\}$ such that
$\delta'(\phi(h,m))-1$ is the number of elements $(h',m') \in J$ such that either
$m' < m-1$ or $m' = m-1$ and $h'\leq h$. Then $\delta'(k) \leq k$ for all
$k \in \{1,...,t\}$, and if $(h,m)$ and $(h,m-1)$ both belong to $J$ then
$\delta'(\phi(h,m)) = \phi(h,m-1) +1$. Also $k_1 < \delta'(k_2)$ if and only if
$k_2 > \delta(k_1)$, and Lemma \ref{5.3} tells us that if $i \in \Delta^{k_1}$
and $j \in \Delta^{k_2}$ then $\beta.(e_i - e_j) \geq |\!| \beta |\!|^2$ if and only
if $k_1 < \delta'(k_2)$. The dual version of Proposition \ref{5.6} tells us that
if $1 \leq k \leq t$ then $E_k/E_{\delta'(k)-1}$ is a direct sum of stable
bundles all of the same slope, and using Remark \ref{4.7} we obtain the 
following dual version of Proposition \ref{5.11}.
\end{rem}

\begin{prop} \label{5.13} 
Let $\beta$ be as in Proposition \ref{5.1} and let $E$ be a semistable bundle
representing an element of $\Sigma_{\beta,l}$ with $\delta$-filtration
$$0=E_0 \subset E_1 \subset ... \subset E_t = E$$
such that if $1 \leq k \leq t$ then $E_k/E_{k-1} \cong \bigoplus_{i=1}^q \CC^{m_i^k} \otimes D_i$
as in Proposition \ref{5.6}. Suppose also that $k_1 \in \{1,...,t\}$ is
such that
$$\beta.e_j >0 \mbox{ whenever } j \in \Delta^{k_2} \mbox{ with } k_2 < \delta'(k_1).$$
Then whenever $F$ is a subbundle of $E$ with $\mbox{slope}(F)= \mbox{slope}(E)$
and such that $E_{{k_1}-1} \subset F \subseteq E_{k_1}$, the subbundle $E_{\delta'(k_1)-1}$
is the minimal subbundle of $F$ such that $F/E_{\delta'({k_1})-1}$
 is a direct sum of stable bundles all with the same slope as $F$.
\end{prop}

\begin{rem} \label{5.14} It follows from the definition of $\Delta^1,...,\Delta^t$
(Definition \ref{5.5}) that if $j_1 \in \Delta^{k_1}$ and $j_2 \in \Delta^{k_2}$ then
$\beta.e_{j_1} < \beta. e_{j_2}$ if and only if $k_1>k_2$, so we can choose $k_-$
and $k_+$ such that $\beta.e_j < 0$ (respectively $\beta.e_j > 0$) if and only if
$j \in \Delta^k$ with $k>k_-$ (respectively $k<k_+$).  Then $k_- = k_+$ or
$k_- = k_+ -1$, depending on whether there exists $j$ with $\beta.e_j = 0$. 
Propositions \ref{5.11} and \ref{5.13} tell us that if $E$ is a semistable bundle
representing an element of $\Sigma_\beta$ with filtration
$$0=E_0 \subset E_1 \subset ... \subset E_t = E$$
as in Proposition \ref{5.6}, then
$$0 \subset E_{\delta(j)+1}/E_{j-1} \subset E_{\delta(\delta(j)+1)+1}/E_{j-1} \subset ... \subset E/E_{j-1}$$
is the maximal Jordan--H\"{o}lder filtration of $E/E_{j-1}$ if $j\geq \delta'(k_-) $, and
$$0 \subset ... \subset E_{\delta'(\delta'(j)-1)-1} \subset E_{\delta'(j)-1} \subset E_{j}$$
is the minimal Jordan--H\"{o}lder filtration of $E_{j}$ if $j \leq \delta(k_+) $.
Note also that $\delta'(k_-) \leq k_- \leq k_+ \leq \delta(k_+)$, so there are values of $j$
satisfying both $j\geq \delta'(k_-) $ and $j \leq \delta(k_+) $.
\end{rem}

There is a converse to Propositions \ref{5.11} and \ref{5.13}.

\begin{prop} \label{conv46} Let 
$\beta$ be as in Proposition \ref{5.1} and let $E$ be a semistable bundle
with $\delta$-filtration
$$0=E_0 \subset E_1 \subset ... \subset E_t = E$$
such that if $1 \leq k \leq t$ then $E_k/E_{k-1} \cong \bigoplus_{i=1}^q \CC^{m_i^k} \otimes D_i$
as in Proposition \ref{5.6}. Suppose that every subbundle $F$ of $E$ with the same slope as
$E$ satisfies the following two properties:

(i) if $E_{{k_1}-1} \subseteq F \subset E_{k_1}$ for some 
$k_1 \in \{1,...,t\}$ 
such that
$\beta.e_j <0$ whenever $ j \in \Delta^{k_2}$ with $ k_2 > \delta(k_1)$,
then the subquotient $E_{\delta(k_1)/F}$
is the maximal subbundle of $E/F$ which 
 is a direct sum of stable bundles all with the same slope as $E/F$;

(ii) if $E_{{k_1}-1} \subset F \subseteq E_{k_1}$ for some 
$k_1 \in \{1,...,t\}$ 
such that
$\beta.e_j >0$  whenever $ j \in \Delta^{k_2}$  with  $k_2 < \delta'(k_1)$, 
then the subbundle $E_{\delta'(k_1)-1}$
is the minimal subbundle of $F$ such that $F/E_{\delta'({k_1})-1}$
 is a direct sum of stable bundles all with the same slope as $F$.

\noindent Then $E$ represents an element of the stratum $\Sigma_{\beta,l}$.
\end{prop}

\noindent{\bf Proof}: Suppose for a contradiction that $E$ does not represent an
element of $\Sigma_{\beta,l}$. Then (cf. \cite{K7}) after applying a change of coordinates
to $\CC^{m_i^k}$  for $1 \leq i \leq q$ and $1 \leq k \leq t$, we can assume that $\beta$
is not equal to the closest point to 0 in the convex hull of $\{e_i - e_j:(i,j) \in S_0\}$ where
$S_0$ is as in the proof of Proposition 8.4. Moreover
$$S_0 \subseteq \{(i,j): \beta.(e_i - e_j) \geq |\!|\beta |\!|^2 \}.$$
Thus $\beta$ does not lie in the convex hull of
$$\{ e_i - e_j:(i,j) \in S_0 \cap S \}$$
where
$S = \{(i,j):  \beta.(e_i - e_j) = |\!|\beta |\!|^2 \}.$
From Lemma \ref{5.3} and Remark \ref{5.4} we know that if $i \in \Delta_{h,m}$ and
$j \in \Delta_{h',m'}$ then $ \beta.(e_i - e_j) \geq |\!|\beta |\!|^2 $ if and only if
$m' \geq m+2$ or $m' = m+1$ and $h' \geq h$, and this happens if and only if
$\phi (h',m') > \delta(\phi(h,m))$, while  $\beta.(e_i - e_j) = |\!|\beta |\!|^2 $
if and only if $m' = m+1$ and $h' = h$. By Remark \ref{5.14} the hypothesis (i)
on subbundles $F$ of $E$ tells us that 
if $E_{{k_1}-1} \subseteq F \subset E_{k_1}$ where
$k_1 \geq \delta(k_-)$, 
then the subquotient $E_{\delta(k_1)/F}$
is the maximal subbundle of $E/F$ which 
 is a direct sum of stable bundles all with the same slope as $E/F$.
This implies that if $(h,m)$ and $(h,m+1)$ both lie in $J$ and $\phi(h,m) \geq k_-$
then every pair $(i,j)$ with $i \in \Delta_{h,m}$ and $j \in \Delta_{h,m+1}$ lies in
$S_0$. Similarly the hypothesis (ii) tells us that if $(h,m)$ and $(h,m-1)$ both lie
in $J$ and $\phi(h,m) \leq k_+$ 
then every pair $(i,j)$ with $i \in \Delta_{h,m}$ and $j \in \Delta_{h,m-1}$ lies in
$S_0$. Since $k_- \leq k_+$ this means that
$$S_0 \cap S = S.$$
This contradicts the fact that $\beta$ is the closest point to 0 in the convex hull
of  $\{ e_i - e_j:(i,j) \in S \}$ but 
does not lie in the convex hull of
$\{ e_i - e_j:(i,j) \in S_0 \cap S \}$, and thus completes the proof.

\begin{rem} \label{new119} Suppose that $\beta$ corresponds to a partition
$\{ \Delta_{h,m}:(h,m) \in J \}$
of $\{1,...,M \}$, indexed by 
$$J = \{(h,m) \in \ZZ \times \ZZ: 1 \leq h \leq L, l_1(h) \leq m \leq l_2(h) \}$$
where $l_1$ and $l_2:\{1,...,L\} \to \ZZ$
satisfy $l_1(h) \leq l_2(h)$ for all $h \in \{1,...,L\}$, as in Proposition \ref{5.1}.
Let $E$ be a semistable bundle representing an element of $\Sigma_{\beta}$
with $\delta$-filtration
$$0=E_0 \subset E_1 \subset ... \subset E_t = E$$
as in Proposition \ref{5.6}. If
$$l_1(h_1) - \epsilon_1(h_1) < l_1(h_2) - \epsilon(h_2) + 1$$
for all $h_1,h_2 \in \{1,..,L\}$, or equivalently if
$$\phi(h,l_1(h)) \leq \delta(1) + 1$$
for all $h \in \{1,...,L\}$, then the proof of  Proposition \ref{conv46} shows that
\begin{equation} \label{nmb} 
0 \subset E_{\delta(1) + 1} \subset E_{\delta(\delta(1) + 1) + 1} \subset ... \subset E
\end{equation}
is the maximal Jordan--H\"{o}lder filtration of $E$. Thus for such a $\beta$ the 
stratum $\Sigma_{\beta,l}$ is contained in the subset $S^{maxJH}_{[{\bf n}_\beta,
{\bf m}_\beta ]}$ of $\mathcal{C}^{ss}$ defined at Definition \ref{def81}, where
${\bf n}_\beta$ and ${\bf m}_\beta$ are determined by the filtration (\ref{nmb}).
If, on the other hand, there exists $h_0 \in \{1,...,L\}$ with
$$\phi(h_0, l_1(h_0)) > \delta(1) + 1,$$
then $E_{\delta(1)+1}$ may not be the maximal subbundle of $E$ which is a
direct sum of stable bundles all with the same slope as $E$; there may be a 
subbundle of $E_{\phi(h_0,l_1(h_0))}/ E_{\phi(h_0,l_1(h_0)-1)}$ which provides
a trivial extension of $E_{\delta(1)+1}$ by a
direct sum of stable bundles all with the same slope as $E$
(see Example \ref{rexlast} below). However, even in this case a careful
analysis of the proof of Proposition \ref{4.4} reveals that it can be modified to
show that the intersection  $\Sigma_\beta \cap S^{maxJH}_{[{\bf n}, {\bf m}]}$ 
is a locally closed complex submanifold of $\mathcal{C}^{ss}$ of finite codimension
for each $\beta \in \Gamma$ and $[{\bf n,m}] \in \mathcal{I}^{ss}$.
\end{rem}

\begin{df} \label{pivot} We shall call a filtration
$0 \subset P_1 \subset ... \subset P_\tau \subset E$ of a semistable bundle
$E$ a {\em pivot} if each subbundle $P_j$ has the same slope as $E$
and $P_1$ is the minimal subbundle of $P_\tau$ such that $P_\tau /P_1$ is 
a direct sum of stable bundles of the same slope, while $P_\tau /P_1$ is the
maximal subbundle of $E/P_1$ which is a direct sum of stable bundles
of the same slope. Any pivot determines a filtration
$$ 0 \subseteq ... \subseteq P_\tau^{-2}=P^{-1}_1 \subseteq 
...\subseteq P^{-1}_{\tau -1} \subseteq P^{-1}_\tau = P_1^0 = P_1 \subseteq 
... $$
$$ ... \subseteq P_\tau = P_\tau^0 =P_1^{+1} \subseteq P_2^{+1} \subseteq ...
\subseteq P_\tau^{+1}  = P_1^{+2} \subseteq ... \subseteq E $$
of $E$ where
$$0 \subset P_j^{+1}/P_j \subset P^{+2}_j/P_j \subset ... \subset E/P_j$$
is the maximal Jordan--H\"{o}lder filtration of $E/P_j$, and
$$0 \subset ... \subset P^{- 2}_j \subset P^{-1}_j \subset P_j$$
is the minimal Jordan--H\"{o}lder filtration of $P_j$.
 A filtration 
$0 = E_0 \subset E_1 \subset ... \subset E_t = E$
of this form (once repetitions have been omitted) for some pivot
$0 \subset P_1 \subset ... \subset P_\tau \subset E$ will be
called a {\em pivotal filtration}. 
It will be called a {\em strongly pivotal filtration} if every subbundle
$F$ of $E$ such that
$$P^m_{j-1} \subseteq F \subset P^m_j$$
for some $m \geq 0$ has
$$0 \subset P^m_j/F \subset P^{m+1}_j/F \subset ... \subset E/F$$
as the maximal Jordan--H\"{o}lder filtration of $E/F$, and
every subbundle
$F$ of $E$ such that
$$P^m_{j} \subset F \subseteq P^m_{j+1}$$
for some $m \leq 0$ has
$$0 \subset ... \subset P^{m-1}_j \subset P^{m}_j  \subset F$$
as its minimal Jordan--H\"{o}lder filtration.

Note that a pivotal
filtration is a $\delta$-filtration where 
$\delta(k_1)$  is the number of $k_2 \in \{1,...,t\}$ 
for which it is not the case that $E_{k_1} = P_{j_1}^{m_1}$ and
$E_{k_2} = P_{j_2}^{m_2}$ with $m_1 \geq m_2$
or $m_1 = m_2 - 1$ and $j_1 \geq j_2$; if the associated
$\delta$-filtration is balanced then we will call the 
pivotal filtration balanced.
\end{df}

\begin{thm}
\label{pivo} Let $\beta$ be as in Proposition \ref{5.1} and let $E$ be a semistable bundle
representing an element of $\Sigma_{\beta,l}$ with $\delta$-filtration
\begin{equation} \label{piveqn} 0=E_0 \subset E_1 \subset ... \subset E_t = E \end{equation}
such that if $1 \leq k \leq t$ then $E_k/E_{k-1} \cong \bigoplus_{i=1}^q \CC^{m_i^k} \otimes D_i$
as in Proposition \ref{5.6}.  Then (\ref{piveqn}) is a balanced strongly pivotal $\delta$-filtration
with pivot
$$0 \subset E_{k_1-1} \subset E_{k_1} \subset ...
\subset E_{\delta(k_1)} \subset E$$
for any $k_1$ satisfying $\delta'(k_-) \leq k_1 \leq k_+$.
\end{thm}
\noindent{\bf Proof}: This is an immediate consequence of Propositions
\ref{5.11} and \ref{5.13} as in Remark \ref{5.14}.

\section{Refinements of the Yang-Mills stratification}

\renorm

We thus have three refinements of the Harder-Narasimhan filtration of a
holomorphic bundle $E$ over $\Sigma$: the maximal Jordan--H\"{o}lder
filtration, the minimal Jordan--H\"{o}lder filtration and the balanced
$\delta$-filtration of maximal triviality obtained by applying Proposition
\ref{5.6} to the subquotients of the Harder-Narasimhan filtration.
Associated to these we have refinements of the Yang-Mills
stratification of
$\mathcal{C}$, each of which is a 
stratification of $\mathcal{C}$  by locally closed complex submanifolds of finite
codimension, and has the set $\mathcal{C}^s$ of stable holomorphic structures
on $\mathcal{E}$ as its open stratum. The first of these refined stratifications is the stratification
\begin{equation} \label{smax} \{ S^{maxJH}_{[{\bf d,n,m}]}: [{\bf d,n,m}] \in \mathcal{I} \}
\end{equation}
of $\mathcal{C}$ defined in Definition 8.1, and another is the stratification
\begin{equation} \label{smin} \{ S^{minJH}_{[{\bf d,n,m}]}: [{\bf d,n,m}] \in \mathcal{I} \}
\end{equation}
defined dually using minimal Jordan--H\"{o}lder filtrations as
in Remark \ref{4.7} and Definition \ref{4.8}. A third refinement is the
stratification
obtained by applying the 
stratification $\{ \Sigma_\gamma: \gamma \in
\Gamma \}$ of $\mathcal{C}^{ss}$, whose indexing set was determined in $\S$5
and whose strata were described in terms of balanced $\delta$-filtrations
in $\S$6, to the $\mathcal{C}(n',d')^{ss}$ which appear inductively in the description 
of the Yang--Mills stratification.

\begin{example} \label{rexlast} Recall from Remark \ref{4.7} that the maximal Jordan--H\"{o}lder
filtration of the direct sum $E \oplus F$ of two semistable bundles of the
same slope is the direct sum of the maximal Jordan--H\"{o}lder filtrations
of $E$ and $F$ with the shorter one extended trivially at the top,
while the minimal Jordan--H\"{o}lder
filtration of $E \oplus F$ is the direct sum of their minimal Jordan--H\"{o}lder filtrations
with the shorter one extended trivially at the bottom. Suppose now that $E$
and $F$ have balanced $\delta$-filtrations of maximal triviality given by
$$0 \subset E_{l_1(1)} \subset E_{l_1(1)+1} \subset ... \subset E_{l_2(1)}=E$$
and
$$0 \subset F_{l_1(2)} \subset F_{l_1(2)+1} \subset ... \subset F_{l_2(2)}=F$$
where the indices $l_1(1),...,l_2(1) \in \ZZ$
and $l_1(2),...,l_2(2) \in \ZZ$ have been chosen so that
$$\epsilon(1) = \sum_{m=l_1(1)}^{l_2(1)} m \mbox{ rank}(E_m/E_{m-1})$$
and
$$\epsilon(2) = \sum_{m=l_1(2)}^{l_2(2)} m \mbox{ rank}(F_m/F_{m-1})$$
lie in the interval $[ - 1/2, 1/2)$. To simplify the notation let us assume that
$l_1(1) \leq l_1(2) \leq l_2(2) \leq l_2(1)$. 
If $\epsilon(1) > \epsilon(2)$ then
$E \oplus F$ has a balanced $\delta$-filtration given by
$$0 \subset E_{l_1(1)} \oplus 0 \subset ... \subset E_{l_1(2)} \oplus 0
\subset E_{l_1(2)} \oplus F_{l_1(2)} \subset ...$$
$$... \subset E_{m-1} \oplus F_{m-1}
\subset E_m \oplus F_{m-1} \subset E_m \oplus F_m \subset ...$$
$$... \subset E_{l_2(2)} \oplus F_{l_2(2)} \subset 
E_{l_2(2)+1} \oplus F_{l_2(2)} \subset  ... \subset 
E_{l_2(1)} \oplus F_{l_2(2)} = E \oplus F.$$
If we assume that
$E_i/E_{i-1}$ and $F_j/F_{j-1}$ are stable for $l_1(1) \leq i \leq l_2(1)$
and $l_1(2) \leq j \leq l_2(2)$, then
this filtration has no proper refinements with subquotients of the same slope as
$E \oplus F$, so by Corollary \ref{nice} it is a balanced $\delta$-filtration of
$E \oplus F$ with maximal triviality (and in fact it is not hard to check that this
is still true without the simplifying assumption).
If $\epsilon(2) > \epsilon(1)$ then we replace the filtration above with the balanced
$\delta$-filtration 
$$0 \subset E_{l_1(1)} \oplus 0 \subset ... \subset E_{l_1(2)-1} \oplus 0
\subset E_{l_1(2)-1} \oplus F_{l_1(2)} \subset ... $$
$$...\subset E_{m-1} \oplus F_{m-1}
\subset E_{m-1} \oplus F_{m} \subset E_m \oplus F_m \subset ...$$
$$... \subset E_{l_2(2)} \oplus F_{l_2(2)} \subset 
E_{l_2(2)+1} \oplus F_{l_2(2)} \subset  ... \subset 
E_{l_2(1)} \oplus F_{l_2(2)} = E \oplus F.$$
Thus we see that the the maximal Jordan--H\"{o}lder
filtration, the minimal Jordan--H\"{o}lder filtration and the balanced
$\delta$-filtration of maximal triviality of a bundle $E$ can all be
different from one another, and that none of them is necessarily a
refinement of the other two. Nonetheless, the concepts of
maximal Jordan--H\"{o}lder filtration, minimal Jordan--H\"{o}lder
filtration and balanced $\delta$-filtration of maximal triviality
on a bundle $E$ are related 
 by Theorem \ref{pivo} via the
notion of a pivotal filtration 
(see also Remark \ref{5.6b}, Propositions \ref{5.11}, \ref{5.13} and
\ref{conv46} and Remark \ref{new119}). 
\end{example}

\bibliographystyle{amsalpha}

\end{document}